\documentclass[11pt]{article}
\usepackage{anysize}
\marginsize{2.0cm}{2.0cm}{2.0 cm}{2.0 cm}
\usepackage{setspace}

\usepackage{latexsym,amsfonts,amsmath}

\usepackage{graphicx}
\usepackage{color}
\usepackage{epstopdf}
\usepackage[titletoc]{appendix}

\def\NN{\mathbb{N}}

\def\RR{\mathbb{R}}

\def\II{\mathbb{I}}

\newcommand{\defi}{:=}

\newtheorem{condition}{Condition}[section]{\bfseries}{\itshape}

{\bfseries}{\itshape}

{\bfseries}{\itshape}

\newtheorem{theorem}{Theorem}[section]{\bfseries}{\itshape}

\newtheorem{corollary}{Corollary}[section]{\bfseries}{\itshape}

\newtheorem{proposition}{Proposition}[section]{\bfseries}{\itshape}

{\bfseries}{\itshape}

\newtheorem{lemma}{Lemma}[section]{\bfseries}{\itshape}

\newtheorem{remark}{Remark}[section]{\bfseries}{\itshape}

\newtheorem{definition}{Definition}[section]{\bfseries}{\itshape}

{\bfseries}{\itshape}

\begin{document}

\author{
Alexey Piunovskiy\thanks{Corresponding author.} \\
Department of Mathematical Sciences, University of Liverpool, L69 7ZL, UK.\\ \texttt{piunov@liv.ac.uk}\\ \ \\
Yi Zhang\\
Department of Mathematical Sciences, University of Liverpool, L69 7ZL, UK.\\
\texttt{yi.zhang@liv.ac.uk}
}
\title{Aggregated occupation measures and linear programming approach to constrained impulse control problems\thanks{Declarations of interest: none.}}
\date{}

\maketitle

\begin{abstract}
For a constrained optimal impulse control problem of an abstract dynamical system, we introduce the occupation measures along with  aggregated occupation measures and present two associated linear programs. We prove that the two linear programs are equivalent under appropriate conditions, and each linear program gives rise to an optimal strategy in the original impulse control problem.
\end{abstract}
\begin{tabbing}
\small \hspace*{\parindent}  \= {\bf Keywords:}
Dynamical System, Optimal Control, Impulse Control, Total Cost, Constraints,\\ Linear Programming.\\
\> {\bf AMS 2000 subject classification:} \= Primary 49N25; Secondary 90C40.
\end{tabbing}

\section{Introduction}

Impulse control of dynamical systems attracts attention of many researchers. The underlying system can be described in terms of ordinary differential equations, see \cite{b9,Barles:1985,Blaquiere:1985,b13,b11,b8,Liu:1998,Miller:2003}, or by a fixed flow in an Euclidean space or in an abstract Borel space, see \cite{DHP,PiunovskiySasha:2018}. An impulse or an intervention means an instantaneous change of the state of the system. In most of the aforementioned works, the target was to optimize a single objective functional, typically having the shape of the integral with respect to time of the running cost and the impulse costs. The popular methods of attack to such problems include dynamic programming, see \cite{b9,Barles:1985,DHP,PiunovskiySasha:2018}, and Pontryagin maximum principle, see \cite{Blaquiere:1985,b8,Miller:2003}. When the total number of impulses is fixed over a finite horizon, the impulse control problem can be treated as a parameter optimization problem, see \cite{b11,Liu:1998}.

In this paper, we consider an impulse control problem of a dynamical system over an infinite horizon with multiple objectives. For optimal control problems with functional constraints, dynamic programming is not always convenient, and the so called convex analytic approach, also known as the linear programming approach, proved to be effective, e.g., for Markov decision processes, see \cite{Dufour:2012,las96,Hernandez-Lerma:1999},  and for deterministic optimal control problems without impulses, see \cite{Gaitsgory:2009,HernandezHernandez:1999,Lasserre:2008}. In a nutshell, this approach, if justified, reduces the original optimal control problem to a linear program in the space of so called occupation measures with the same (optimal) value, and one can retrieve an optimal control strategy for the original problem from the optimal solution to the induced linear program.

For a deterministic impulse control problem over a finite horizon, a linear program formulation was presented in \cite{b13}, from which, as the primitive goal of that paper, the authors established a numerical method for solving approximately the original problem. For this reason, \cite{b13} dealt with an unconstrained problem for a specific model with polynomial initial data, and did not show that the formulated linear program was equivalent to the original impulse control problem. 
Another, slightly different linear programming approach appeared in \cite{b13p,b13pp}, where the equivalence between the linear program and the original problem was briefly discussed. In the aforementioned works, the flow in an Euclidean space came from an ordinary differential equation, whereas in the present paper the flow is arbitrary enough and lives in a Borel space.
A different linear program formulation  was presented in \cite{ZYnew}, which was shown to be equivalent to the original impulse control problem.

In this article, we start with recapitulating briefly the linear programming approach developed in \cite{ZYnew}, which was in the space of occupation measures, see (\ref{SashaLp02}) and (\ref{e106}). Then we introduce the second linear program, which is in the space of so called aggregated occupation measures and is connected to the specific linear programs  described in \cite{b13,b13p,b13pp}. As the term suggests, aggregated measures arise from suitably aggregating the occupation measures, see (\ref{e10}), (\ref{e11}). The main difference and advantage of the aggregated occupation measures are in the reduction of the dimensionality: see Remark \ref{rem991}.
Our main contributions lie in that we prove the equivalence between the  mentioned above  linear programs, see Corollary \ref{corol1}, and show that the ``induced'' strategy from either one solves the original impulse control problem.

In simple cases (see Section \ref{sec7}), the second linear program, after the suitable change of measures, can be transformed to the linear programs obtained in \cite{b13,b13p,b13pp}. The novelty of the present article is in the following:
\begin{itemize}
\item the dynamical system is described by a flow in an arbitrary Borel space, rather than by an ordinary differential equation in an Euclidean space;
\item the optimal solution must satisfy a number of functional constraints which were absent in the cited literature;
\item under suitable conditions, we rigorously prove that the optimal values of the original impulse control problem and of the introduced linear programs coincide, i.e., there is no ``relaxation gap'';
\item we show how to retrieve the optimal control strategy from the solutions to the  associated linear programs.
\end{itemize}

The rest of this article is organized as follows. The problem statement is described in Section \ref{sec2}. In Sections \ref{sec3} and \ref{sec33}, we formulate the preliminary observations and the main results  correspondingly. In Section \ref{sec7}, we present an example and compare our approach with works \cite{b13,b13p,b13pp}. The proofs of the main theorems are given in
 Sections \ref{sec5} and \ref{sec6}. Some auxiliary lemmas are presented and proved in the Appendix.

Throughout this paper, we use the following notations: $\NN:=\{1,2,\ldots\}$, $\bar{\mathbb{R}}_+^0:=[0,\infty]$, $\mathbb{R}_+^0:=[0,\infty)$, $\mathbb{R}_+:=(0,\infty)$.  The term ``measure" will always refer to a countably additive $\bar\RR^0_+$-valued set function, equal to zero on the empty set. Consider two $\sigma$-finite measures $\eta_1$ and $\eta_2$ on a common measurable space $(\Omega,{\cal F})$ such that $\eta_1\ge \eta_2$ set-wise. Then there exists a measurable decomposition $\{\Omega_n\}_{n=1}^\infty$ of $\Omega$ such that ${\eta}_1( \Omega_n)<\infty$ and ${\eta}_2(\Omega_n)<\infty.$ The difference between these two measures is defined by $(\eta_1-\eta_2)(d\omega):=\sum_{n=1}^\infty (\eta_1 (d\omega \cap \Omega_n)-\eta_2(d\omega \cap \Omega_n))$. ${\cal P}(E)$ is the space of all probability measures on a measurable space $(E,{\cal B}(E))$. On the time axis $\RR^0_+$  the expression ``for almost all $u$'' is understood with respect to the Lebesgue measure.
By default, the $\sigma$-algebra on $\RR^0_+$ is just the Borel one.
If $(E,{\cal B}(E))$ is a measurable space then, for $Y\in{\cal B}(E)$, ${\cal B}(Y)\defi\{X\cap Y,~X\in{\cal B}(E)\}$ is the restriction of the $\sigma$-algebra ${\cal B}(E)$. Integrals on a measure space $(E,{\cal B}(E),\mu)$ are denoted as $\int_E h(e)d\mu(e)$ or as $\int_E h(e) \mu(de)$.  If $b=\infty$ then the Lebesgue integrals $\displaystyle \int_{[a,b]}f(u)du$ are taken over the open interval $(a,\infty)$. Expressions like ``positive, negative, increasing, decreasing'' are understood in the non-strict sense, like ``nonnegative'' etc. For $I\subset \RR$, $\tau\in\RR$, $\tau+I\defi\{\tau+x:~x\in I\}$ is the shifted set. $\II\{\cdot\}$ is the indicator function; $\delta_y(dx)$ is the Dirac measure at the point $y$. For $b,c\in[-\infty,+\infty]$,   $b^+:=\max\{b,0\}$, $b^-:=-\min\{b,0\}$, $b\wedge c:=\min\{b,c\}$, $b\vee c:=\max\{b,c\}$.

\section{Problem Statement}\label{sec2}

We will deal with a control model defined through the following elements.
\begin{itemize}
\item ${\bf X}$ is the state space, which is a topological Borel space.
\item $\phi(\cdot , \cdot):~{\bf X}\times\RR^0_+\to{\bf X}$ is the measurable flow possessing the semigroup property $\phi( x,t+s)=\phi(\phi( x,s),t)$ for all $ x\in{\bf X}$ and $(t,s)\in(\RR^0_+)^2$; $\phi( x,0)= x$ for all $ x\in{\bf X}$. Between the impulses, the state changes according to the flow.
\item $\bf A$ is the action space, again a topological Borel space with a compatible metric $\rho_A$.
\item $ l(\cdot , \cdot):~{\bf X}\times{\bf A}\to{\bf X}$ is the mapping describing the new state after the corresponding action/impulse is applied.
\item For each $j=0,1,\dots,J,$ where and below $J$ is a fixed natural number, $C^g_j(\cdot):~{\bf X}\to\RR_+^0$ is a (gradual) cost rate.
\item For each $j=0,1,\dots,J,$ $C^I_j(\cdot\,, \cdot):~{\bf X}\times{\bf A}\to\RR_+^0$ is a cost function associated with the actions/impulses applied in the corresponding states.
\end{itemize}
All the mappings $\phi, l,\{C^g_j\}_{j=0}^J$ and $\{C^I_j\}_{j=0}^J$ are assumed to be measurable. The initial state $ x_0\in{\bf X}$ is fixed.

We assume that the states $x\in{\bf X}$ have the form $x=(\tilde x,t)$, where $t\in\RR^0_+$ equals time elapsed since the most recent impulse, and $\tilde x\in\tilde{\bf X}$, an arbitrary Borel space with a compatible metric $\tilde\rho$.
In this connection,
$$\phi(x,u)=\phi((\tilde x,t),u)\defi(\tilde\phi(\tilde x,u),t+u),$$
where $\tilde\phi(\cdot , \cdot):~\tilde{\bf X}\times\RR^0_+\to\tilde{\bf X}$ is the measurable flow in $\tilde{\bf X}$ possessing the semigroup property. Similarly, $l(x,a)=(\tilde l(x,a),0)$, where $ \tilde l(\cdot, \cdot):~{\bf X}\times{\bf A}\to\tilde{\bf X}$ is a measurable mapping:
after each impulse, the $t$-component goes down to zero. Any initial state is in the form $x_0=(\tilde x_0,0)$ and thus has the time component zero. The mappings $\tilde\phi$ and $\tilde l$ are assumed to be measurable.

\begin{remark}\label{remark2}
If the original state space is just $\tilde{\bf X}$, then it is always possible to extend it by including the component $t$.
\end{remark}

We exclude from the consideration all the points from  $\tilde{\bf X}\times\RR^0_+$ which cannot appear in the dynamical system generated by the  flow $\tilde\phi$,  so that
$${\bf X}:=\{(\tilde y,t)\in\tilde{\bf X}\times\RR^0_+:~~\tilde y=\tilde\phi(\tilde x,t)~\mbox{ for some } \tilde x\in\tilde{\bf X}\}.$$
In $\RR^0_+$, the standard Euclidean topology is fixed.
The product space $\tilde{\bf X} \times \mathbb{R}^0_+$ is equipped with the product topology, which is metrizable (see \cite[\S2.14]{b4}). The topology on $\textbf{X}$ is the restriction of the product topology on $\tilde{\bf X} \times \mathbb{R}^0_+$  on it. We endow $\textbf{X}$ with its Borel $\sigma$-algebra, which is the restriction of the Borel $\sigma$-algebra ${\cal B}(\tilde{\bf X} \times \mathbb{R}^0_+)$ on $\textbf{X}$, see \cite[Lem.7.4]{Bertsekas:1978}. Since $\textbf{X}$ is a projection of the graph of the measurable mapping $\tilde{\phi}$, it is not immediately obvious whether $\textbf{X}$ is a Borel subset of $\textbf{X}\times \mathbb{R}_+^0$. In this and the next section, we assume that ${\bf X}$ is a Borel space. Sufficient conditions will be imposed later to guarantee this is indeed the case (see Lemma \ref{la5}).

Let ${\bf X}_\Delta:={\bf X}\cup\{\Delta\}$, where $\Delta$ is an isolated artificial point describing the case that the controlled process is over and no future costs will appear. The dynamics (trajectory) of the system can be represented as one of the following sequences
\begin{eqnarray}
&&x_0\to (\theta_1,a_1)\to x_1\to (\theta_2,a_2)\to \ldots;~~~~ \theta_i<+\infty \mbox{ for all } i\in\NN, \nonumber\\
\mbox{or}&&\label{e1}\\
&&x_0\to (\theta_1,a_1)\to\ldots\to x_n\to (+\infty,a_{n+1})\to \Delta\to (\theta_{n+2},a_{n+2}) \to \Delta \to \ldots,\nonumber
\end{eqnarray}
where $x_0\in{\bf X}$ is the initial state of the controlled process and $\theta_i<+\infty$ for all $i=1,2,\ldots, n$. For the state $x_{i-1}\in{\bf X}$, $i\in\NN$, the pair $(\theta_i,a_i)\in\bar\RR^0_+\times{\bf A}=:{\bf B}$ is the control at the step $i$: after $\theta_i$ time units, the impulsive action $a_i$ will be applied leading to the new state
\begin{equation}\label{e1p}
x_i=\left\{\begin{array}{ll}
l(\phi(x_{i-1},\theta_i),a_i), & \mbox{ if } \theta_i<+\infty; \\
\Delta, & \mbox{ if } \theta_i=+\infty.
\end{array}\right.
\end{equation}
The state $\Delta$ will appear forever, after it appeared for the first time, i.e., it is absorbing.

\begin{remark}\label{rem0}
We underline that all the realized points $x_i$, $i=1,2,\ldots$, provided that they are not equal to $\Delta$, have the form $(\tilde x_i,0)$. For technical needs, unless stated otherwise, we allow $x_0$ to be an arbitrary point in ${\bf X}$.
\end{remark}

After each impulsive action, if $\theta_1,\theta_2,\ldots,\theta_{i-1}<+\infty$, the decision maker has in hand the complete information about the history, that is, the sequence
$$x_0, (\theta_1,a_1), x_1,\ldots,  (\theta_{i-1},a_{i-1}),x_{i-1}.$$
The selection of the next control $(\theta_i,a_i)$ is based on this information, and we also allow the selection of the pair $(\theta_i,a_i)$ to be randomized. Below, the control $(\theta,a)\in{\bf B}$ is denoted as $b$.

For each $j=0,1,\dots,J,$ the cost accumulated on the coming interval of length $\theta_i$ equals
\begin{equation}\label{e1pprim}
\int_{[0,\theta_i] } C_j^g(\phi(x_{i-1},u))du+\II\{\theta_i<+\infty\} C^I_j(\phi(x_{i-1},\theta_i),a_i),
\end{equation}
the last term being absent if $\theta_i=+\infty$.
The next state $x_i$ is given by formula (\ref{e1p}).

In the space of all the trajectories (\ref{e1})
\begin{eqnarray*}
\Omega&=&\cup_{n=1}^\infty[{\bf X}\times((\RR^0_+\times{\bf A})\times{\bf X})^n
\times(\{+\infty\}\times{\bf A})\times\{\Delta\}\times ((\bar\RR^0_+\times{\bf A})\times\{\Delta\})^\infty]\\
&& \cup [{\bf X}\times((\RR^0_+\times{\bf A})\times{\bf X})^\infty],
\end{eqnarray*}
we fix the natural $\sigma$-algebra $\cal F$.
Finite sequences
$$h_i=(x_0, (\theta_1,a_1), x_1, (\theta_2,a_2),\ldots,x_i)=(x_0,b_1,x_1,b_2,\ldots,x_i)$$
will be called (finite) histories; $i=0,1,2,\ldots$, and the space of all such histories will be denoted as ${\bf H}_i$; ${\cal F}_i:={\cal B}({\bf H}_i)$ is the restriction of ${\cal F}$ to ${\bf H}_i$. Capital letters $X_i,T_i,\Theta_i, A_i,B_i=(\Theta_i,A_i)$ and $H_i$ denote the corresponding functions of $\omega\in\Omega$, i.e., random elements.

\begin{definition}\label{d1}
A control strategy $\pi=\{\pi_i\}_{i=1}^\infty$ is a sequence of stochastic kernels $\pi_i$ on ${\bf B}=\bar\RR^0_+\times{\bf A}$ given ${\bf H}_{i-1}$.
A Markov strategy is defined by stochastic kernels $\{\pi_i(db|x_{i-1})\}_{i=1}^\infty$.
A control strategy is called stationary, and denoted as $\widetilde{\pi}$, if there is a stochastic kernel $\widetilde{\pi}$ on $\bar\RR^0_+\times{\bf A}$ given $\textbf{X}_\Delta$  such that $\pi_i(db|h_{i-1})=\widetilde{\pi}(db|x_i)$ for all $i=1,2,\ldots$. Every measurable mapping $f:~{\bf X}_\Delta\to{\bf B}$ defines a deterministic stationary strategy, which is given by $\pi_i(db|h_{i-1}):=\delta_{f(x_{i-1})}(db)$, and identified with $f$.
\end{definition}

Note that every Markov strategy can be represented as
$$\pi_i(d\theta\times da|x)=p^i_T(d\theta|x)p^i_A(da|x,\theta),$$
where $p^i_T$ and $p^i_A$ are stochastic kernels on $\bar\RR^0_+$ given ${\bf X}_\Delta$ and on $\bf A$ given ${\bf X}_\Delta\times\bar\RR^0_+$, correspondingly: see \cite[Prop.7.27]{Bertsekas:1978}.

For a given initial state $x\in{\bf X}$ and a strategy $\pi$, there is a unique probability measure $P^\pi_{x}(\cdot)$ on $\Omega$ constructed using the Ionescu-Tulcea Theorem, satisfying
for all $i\in\NN$, $\Gamma\in{\cal B}(\bar\RR^0_+\times{\bf A})$, $\Gamma_X\in{\cal B}({\bf X}_\Delta)$,
\begin{eqnarray}
P^\pi_{x}(X_0\in\Gamma_X)&=&\delta_{x}(\Gamma_X)\mbox{ for } \Gamma_X\in{\cal B}({\bf X}_\Delta);\nonumber\\
P^\pi_{x}((\Theta_i,A_i)\in\Gamma|H_{i-1})&=&\pi_i(\Gamma |H_{i-1});\label{e33}\\
P^\pi_{x}(X_i\in\Gamma_X|H_{i-1},(\Theta_i,A_i))
&=& \left\{\begin{array}{ll}
\delta_{l(\phi(X_{i-1},\Theta_i),A_i)}(\Gamma_X), & \mbox{ if } X_{i-1}\in{\bf X},~\Theta_i<+\infty; \\
\delta_\Delta(\Gamma_X) & \mbox{ otherwise.} \end{array} \right. \nonumber
\end{eqnarray}
This is a standard definition of strategic measures in Markov Decision Processes. Let $E^\pi_{x}$ be the corresponding mathematical expectation.

Let us introduce the notation
\begin{eqnarray*}
&&{\cal V}_j(x, \pi) \\
&:=&  E^\pi_{x}\left[\sum_{i=1}^\infty \II\{X_{i-1}\ne\Delta\} \left\{ \int_{[0,\Theta_{i}]}  C^g_j(\phi(X_{i-1},u)) du\right.\right.+ \II\{\Theta_i<+\infty\} \left.\left.\vphantom{\sum_{i=1}^\infty} C^I_j(\phi(X_{i-1},\Theta_i),A_i)\right\}\right]
\end{eqnarray*}
for each strategy $\pi$, $j=0,1,\dots,J$ and initial state $x\in{\bf X}$.

The constrained optimal control problem under study is the following one:
\begin{eqnarray}\label{PZZeqn02}
\mbox{Minimize with respect to } \pi &&{\cal V}_0(x_0, \pi)  \\
\mbox{subject to }&& {\cal V}_j(x_0,\pi)\le d_j,~j=1,2,\dots,J.\nonumber
\end{eqnarray}
Here and below,  $\{d_j\}_{j=1}^J$ are fixed constraint constants and $x_0=(\tilde x_0,0)$ is a fixed initial state, where $\tilde x_0\in\tilde{\bf X}$.

\begin{definition}\label{d2}
A strategy $\pi$ is called feasible if it satisfies all the constraint inequalities in problem (\ref{PZZeqn02}). A feasible strategy $\pi^\ast$ is called optimal if, for all feasible strategies $\pi,$ ${\cal V}_0(x_0,\pi^*)\le  {\cal V}_0(x_0,\pi)$.
\end{definition}

We shall assume that problem (\ref{PZZeqn02}) is consistent.
\begin{condition}\label{ConstrainedPPZcondition05}
There exists some feasible strategy $\pi$ such that ${\cal V}_0(x_0,\pi)<\infty.$
\end{condition}

In what follows, we develop the linear programming approach to problem (\ref{PZZeqn02}).


\section{Preliminary Observations}\label{sec3}

Clearly, the control model presented in Section \ref{sec2}, from the formal viewpoint, is a specific constrained Markov Decision Process \cite{Altman:1999,Dufour:2012,Hernandez-Lerma:1999,Piunovskiy:1997} , which is defined by the following elements.  The state space is
${\bf X}_\Delta:=\textbf{X}\cup \{\Delta\}$,
as before, where the state $\Delta\notin \textbf{X}$ is an isolated point and $\bf X$ is assumed to be a Borel space. The action space is
$\textbf{B}:=\bar\RR^0_+\times{\bf A}$,
which is endowed with the product topology and the corresponding Borel $\sigma$-algebra. The transition kernel is defined by
\begin{eqnarray*}
Q(dy|x,(\theta,a)):=\left\{\begin{array}{ll}
\delta_{l(\phi(x,\theta),a)}(dy), & \mbox{ if } x\ne\Delta,~\theta\ne+\infty;\\ \delta_\Delta(dy) & \mbox{ otherwise}, \end{array}\right..
\end{eqnarray*}
The cost functions are given by
\begin{eqnarray*}
\bar C_j(x,(\theta,a)):= \II\{x\ne\Delta\} \left\{ \int_{[0,\theta]}  C^g_j(\phi(x,u)) du
+ \II\{\theta<+\infty\} C^I_j(\phi(x,\theta),a)\right\},~j=0,1,\dots,J.
\end{eqnarray*}
If $\theta=\infty$, then the above integration is understood over $[0,\infty)$. Below, we omit such remarks.
The initial state $(\tilde x_0,0)\in{\bf X}$
and the constraint constants  $d_j\in\mathbb{R}_+^0,~j=1,2,\dots,J$ are as before.

Let us impose the next set of compactness-continuity conditions.

\begin{condition}\label{ConstrainedPPZcondition01}
\begin{itemize}
\item[(a)] The space $\bf A$ is compact, and $+\infty$ is the one-point compactification of the positive real line $\RR^0_+$.
\item[(b)] The mapping $(x,a)\in \textbf{X}\times \textbf{A}\rightarrow l(x,a)$ is continuous.
\item[(c)] The mapping $(x,\theta)\in \textbf{X}\times \mathbb{R}_+^0\rightarrow \phi(x,\theta)$ is continuous.
\item[(d)] For each $j=0,1,\dots,J,$ the function $(x,a)\in \textbf{X}\times \textbf{A}\rightarrow C_j^I(x, a)$ is lower semicontinuous.
\item[(e)] For each $j=0,1,\dots,J,$ the function $x\in \textbf{X}\rightarrow C_j^g(x)$ is lower semicontinuous.
\end{itemize}
\end{condition}

According to Theorem 1 of \cite{PiunovskiySasha:2018},
under Condition \ref{ConstrainedPPZcondition01}, assuming that $\textbf{X}$ is a Borel space, the function on $\textbf{X}$ defined by
$
\inf_\pi E^\pi_x\left[\sum_{i=0}^\infty\sum_{j=0}^J \bar C_j(X_i,B_{i+1})\right]
$
is lower semicontinuous.

\begin{condition}\label{co31} There exists $\delta>0$ such that $\sum_{j=0}^J C^I_j(x,a)\ge \delta$ for all $(x,a)\in{\bf X}\times{\bf A}$.
\end{condition}

The above condition asserts that each impulse is costly. Below in this section, we assume that Conditions \ref{ConstrainedPPZcondition01} and \ref{co31} are satisfied.

Consider a point $x\in{\bf X}$ such that
\begin{eqnarray}\label{JulyEqn01}
\inf_\pi E^\pi_x\left[\sum_{i=0}^\infty\sum_{j=0}^J \bar C_j(X_i,B_{i+1})\right]=0
\end{eqnarray}
(provided that such a point exists). Then
$E^{f^*}_x\left[\sum_{i=0}^\infty\sum_{j=0}^J \bar C_j(X_i,B_{i+1})\right]=0
$
for the deterministic stationary strategy
$f^*(x)\equiv (\infty,\hat a)$  with the immaterial value of $\hat a\in{\bf A}$ being arbitrarily fixed: for all other values of $B_1\in{\bf B}$, $\sum_{j=0}^J\bar C_j(x,B_1)\ge\delta>0$.

Clearly, the control $(\infty,\hat a)$ is optimal in problem (\ref{PZZeqn02}) at all such states $x\in{\bf X}$, at which (\ref{JulyEqn01}) holds.
Moreover, for all such states $x$, $Q(\{\Delta\}|x,f^*(x))=1$ and $X_1=\Delta$ $P_x^{f^\ast}$-almost surely, so that
\begin{eqnarray*}
0=E^{f^*}_x\left[\sum_{i=0}^\infty\sum_{j=0}^J \bar C_j(X_i,B_{i+1})\right]= \sum_{j=0}^J \bar C_j(x,(\infty,\hat a))=\int_{[0,\infty)}\sum_{j=0}^J C^g_j(\phi(x,u))du,
\end{eqnarray*}
and consequently, for all $j=0,1,\ldots,J$, $C^g_j(\phi(x,u))=0$  for almost all $u\ge 0$.
Conversely, if, at some $x\in\textbf{X}$, for all $j=0,1,\ldots,J$, $C^g_j(\phi(x,u))=0$ for almost all  $u\ge 0$, then (\ref{JulyEqn01}) holds.

Below, let us denote
\begin{eqnarray*}
V:=\left\{x\in{\bf X}:~\inf_\pi E^\pi_x\left[\sum_{i=0}^\infty\sum_{j=0}^J \bar C_j(X_i,B_{i+1})\right]>0\right\}=\left\{ x\in{\bf X}:\int_{[0,\infty)} \sum_{j=0}^J C^g_j(\phi(x,u))du>0\right\}
\end{eqnarray*}
and $V^c:={\bf X}_\Delta\setminus V$. The set $V$, as the preimage of an open set under a lower semicontinuous function, is open in $\textbf{X}$.  The set $V^c$ can be equivalently defined as
\begin{eqnarray*}
V^c:=\{\Delta\}\cup\left\{ x\in{\bf X}:~ C^g_j(\phi(x,u))=0 \mbox{ for all $j=0,1,\ldots,J$, for almost all $u\ge 0$}\right\},
\end{eqnarray*}
and it is absorbing with respect to the flow $\phi$: for each $x\in{\bf X}$, as soon as $\phi(x,u)\in V^c$, $\phi(x,s)\in V^c$ for all $s\ge u$.
The case $V={\bf X}$ and $V^c=\{\Delta\}$ is not excluded.

In view of the previous observations, under Conditions \ref{ConstrainedPPZcondition01} and \ref{co31}, it is sufficient to
consider the class of reasonable strategies $\pi=\{\pi_i\}_{i=1}^\infty$ defined as follows.

\begin{definition}\label{JulyRem01}
 Assume $\textbf{X}$ is a Borel space, and suppose Conditions \ref{ConstrainedPPZcondition01} and \ref{co31} are satisfied. A strategy $\pi=\{\pi_i\}_{i=1}^\infty$ is called reasonable if
$
\pi_{i}(db|x_0,b_1,x_1,\dots,x_{i-1})=\delta_{f^\ast(x_{i-1})}(db)$ for all $x_{i-1}\in V^c,$ and
\begin{eqnarray*}
\pi_i([\tilde\theta^*(\tilde x_{i-1}),\infty)\times{\bf A}|x_0,b_1,x_1,\ldots,x_{i-1})=0,~i=1,2,\ldots.
\end{eqnarray*}
Here, $x_{i-1}=(\tilde x_{i-1},0)$ (see Remark \ref{rem0}) and
\begin{equation}\label{e14pp}
\tilde\theta^*(\tilde x)\defi\inf\{\theta\in\RR^0_+:~\phi((\tilde x,0),\theta)\in V^c\}
\end{equation}
is a function defined for each $\tilde x\in\tilde{\bf X}$.
(As usual, $\inf \emptyset \defi +\infty$.)
\end{definition}

Since the flow $\phi$ is continuous, the function $\tilde\theta^*(\cdot)$ is measurable: see \cite[Lemma 27.1]{b1} or \cite[Prop.1.5, p.154]{b3}.
After we introduce notations
$$\tilde{V}:=\{\tilde x\in\tilde{\bf X}:~(\tilde x,0)\in{V}\}~\mbox{ and } \tilde{V}^c:=\{\tilde x\in\tilde{\bf X}:~(\tilde x,0)\in{V}^c\}=\tilde{\bf X}\setminus\tilde{V},$$
it is clear that, for $\tilde x\in\tilde V$, $\tilde\theta^*(\tilde x)>0$ because the set $V$ is open and the set $V^c$ is closed; in case $\tilde\theta^*(\tilde x)<\infty$, the infimum in (\ref{e14pp}) is attained, and
\begin{eqnarray*}
\tilde\theta^*(\tilde x)=\sup\{t\in\RR^0_+:~\phi((\tilde x,0),t)\in V\}.
\end{eqnarray*}
If $\tilde x\in\tilde V^c$, then $\tilde\theta^*(\tilde x)=0$.

We thus concentrate on selecting actions at the states $x\in V$ and restrict ourselves to the set of reasonable strategies.

A linear programming method was established in \cite{ZYnew} regarding how to select actions at $x\in V$, and it serves the beginning of the analysis in the present paper. For this reason, let us briefly describe it: see (\ref{SashaLp02}), (\ref{e105}) below. The formulation of that linear program is related to the occupation measures $\mu^\pi$ defined as follows:
\begin{equation}\label{e101}
\mu^\pi(\Gamma_1\times\Gamma_2):=E_{x_0}^\pi\left[\sum_{i=0}^\infty \II\{(X_i,B_{i+1})\in\Gamma_1\times\Gamma_2\}\right],~\forall~\Gamma_1\in{\cal B}(\textbf{X}_\Delta),\Gamma_2\in{\cal B}(\bar\RR^0_+\times{\bf A}).
\end{equation}
Under Conditions \ref{ConstrainedPPZcondition05}, \ref{ConstrainedPPZcondition01} and \ref{co31}, for each reasonable $\pi$ as in Definition \ref{JulyRem01},
\begin{eqnarray*}
{\cal V}_j(x_0,\pi)=\int_{V\times \bar\RR^0_+\times\textbf{A}}\bar{C}_j(x,(\theta,a))\mu^\pi(dx\times d\theta\times da).
\end{eqnarray*}
It follows that the restriction on $V\times\bar\RR^0_+\times{\bf A}$ of any occupation measure $\mu^\pi=\mu$ of our interest is concentrated on the measurable subset ${\bf M}\times{\bf A}$, where
\begin{eqnarray}\label{e1101}
{\bf M}:=\{(y,\theta):  y=(\tilde x,0) \mbox{ with } \tilde x \in \tilde V \mbox{~and~} \theta\in [0,\tilde\theta^*(\tilde x))\cup\{\infty\}\}.
\end{eqnarray}
Moreover, there is no need to consider such occupation measures that $\mu^\pi({\bf M}\times{\bf A})=\infty$: the latter means that, with positive probability, actions from  $\RR^0_+\times{\bf A}$ at states from $V$ appear infinitely many times, leading to the infinite value of at least one of the objectives ${\cal V}_j(x_0,\pi)=\int_{V\times \bar\RR^0_+\times\textbf{A}}\bar{C}_j(x,(\theta,a))$\linebreak$\times\mu^\pi(dx\times d\theta\times da)$.

The impulse control problem (\ref{PZZeqn02}) is now equivalent to the following linear program:
\begin{eqnarray}\label{SashaLp02}
\mbox{Minimize}&:&\int_{V \times\bar\RR^0_+\times\textbf{A}}\bar{C}_0(x,(\theta,a)) \mu(dx\times d\theta\times da) \\
&&\mbox{over finite measures $\mu$ on $V\times{\bf B}=V\times\bar\RR^0_+\times{\bf A}$
concentrated on
${\bf M}\times{\bf A}$}\nonumber
\end{eqnarray}
\begin{eqnarray}
 \mbox{subject to}~ \mu(dx\times\bar\RR^0_+\times \textbf{A})=\delta_{x_0}(dx)+\int_{V\times \bar\RR^0_+\times{\bf A}}Q(dx|y,(\theta,a)) \mu(dy\times d\theta\times da)~\mbox{on ${\cal B}(V)$};\label{e105}\\
\int_{V \times\bar\RR^0_+\times\textbf{A}}\bar{C}_j(x,(\theta,a))\mu(dx\times d\theta\times da)\le d_j,~j=1,2,\dots,J. \nonumber
\end{eqnarray}
See Proposition \ref{pr1} for a precise statement of this equivalence.

One can recognize that the form of this linear program is standard for the total cost Markov Decision Processes (see e.g., \cite{Altman:1999,Dufour:2012,Hernandez-Lerma:1999}). For every reasonable strategy $\pi$, the occupation measure $\mu^\pi$ satisfies equality (\ref{e105}).

The next statement comes from Theorem 4.1 of \cite{ZYnew}.
\begin{proposition}\label{pr1}
Suppose the space $\bf X$ is Borel and  Conditions  \ref{ConstrainedPPZcondition05},  \ref{ConstrainedPPZcondition01} and \ref{co31} are satisfied. Then the following assertions hold.
\begin{itemize}
\item[(a)] There exists a solution $\mu^*$ to the program (\ref{SashaLp02}), which gives rise to the optimal (in problem (\ref{PZZeqn02}))  stationary strategy $\widetilde\pi$ coming from the decomposition
\begin{eqnarray*}
\mu^*(dx\times db)=\mu^*(dx\times{\bf B})\times\widetilde\pi(db|x),~x\in V.
\end{eqnarray*}
On the space $V^c$, the optimal strategy is given by $f^*(x)\equiv(\infty,\hat a)$ as before; the value of $\hat a\in{\bf A}$ is immaterial. The minimal value of the program (\ref{SashaLp02}) is finite and coincides with the minimal value of the original problem (\ref{PZZeqn02}).
\item[(b)] If $\pi^*$ is a reasonable strategy, whose occupation measure $\mu^{\pi^*}$ on $V\times \textbf{B}$ is concentrated on ${\bf M}\times{\bf A}$ and solves the linear program (\ref{SashaLp02}), then the strategy $\pi^*$ is optimal in problem (\ref{PZZeqn02}).
\end{itemize}
\end{proposition}

In what follows, we use the notation
${\bf A}_\Box={\bf A}\cup\{\Box\}$,
where $\Box\notin{\bf A}$ is an artificial isolated point.

The target of this article is to pass to the equivalent in some sense  linear program in the space of  measures $\eta$ on $V\times{\bf A}_\Box$. The reason is connected with the form of the objectives ${\cal V}_j(x_0,\pi)$. Since they are linear with respect to the original functions $C^g_j$ and $C^I_j$ on ${V}$ and ${V}\times{\bf A}$ correspondingly, it is desirable to represent them in the form of $\int_{{V}\times{\bf A}_\Box} C_j(y,a)\eta(dy\times da)$, where
\begin{equation}\label{e103}
C_j(y,a):=\left\{\begin{array}{ll}
C^g_j(y), & \mbox{ if } a=\Box; \\ C^I_j(y,a), & \mbox{ if } a\in{\bf A},
\end{array}\right.
\end{equation}
and develop the characteristic equation for the measures $\eta$.

Consider a finite measure $\mu$ in the linear program (\ref{SashaLp02}), which can be written in the form
\begin{equation}\label{e14p}
\mu(dx\times d\theta\times da)=p_T(d\theta|x,a)p_A(da|x)\mu(dx\times\bar\RR^0_+\times{\bf A}),
\end{equation}
where $p_T(\cdot)$ and $p_A(\cdot)$ are stochastic kernels on $\bar\RR^0_+$ and $\bf A$ correspondingly: see \cite[Prop.7.27]{Bertsekas:1978}. The dependence of $p_T$ and $p_A$ on $\mu$ is not explicitly indicated here. Hence, using the Tonelli Theorem (see \cite[Thm.11.28]{b4}), some straightforward calculations imply that
\begin{eqnarray*}
&&\int_{V\times \bar\RR^0_+\times \textbf{A}} \left\{ \int_{[0,\theta]}  C^g_j(\phi(x,u)) du
\right\}\mu(dx\times d\theta\times da)\\
&=&\int_{V\times \bar\RR^0_+\times \textbf{A}} \left\{ \int_{[0,\theta]}  C^g_j(\phi(x,u))\II\{\phi(x,u)\in V\} du
\right\}\mu(dx\times d\theta\times da)\\
&=& \int_{V}\int_{\bf A}\int_{\bar\RR^0_+}\int_{[0,\theta]}  C^g_j(\phi(x,u))\II\{\phi(x,u)\in V\} du~ p_T(d\theta|x,a)p_A(da|x)\mu(dx\times \bar\RR^0_+\times{\bf A})\\
&=& \int_{V}\int_{\bf A}\int_{\RR^0_+} C^g_j(\phi(x,u))\II\{\phi(x,u)\in V\}  p_T([u,\infty]|x,a)~du~p_A(da|x)\mu(dx\times \bar\RR^0_+\times{\bf A}),
\end{eqnarray*}
where the first equality holds because $C^g_j(\phi(x,u))=0$ for each $x\in V$, for almost all $u\in\RR^0_+$ such that $\phi(x,u)\in V^c$,  for all $j=0,1,\ldots, J$. To put it differently, $C^g_j(\phi(x,u))=0$ for almost all $u\ge\tilde\theta^*(\tilde x)$ for all $x=(\tilde x,0)$.

After we introduce the following measure on $V$
\begin{eqnarray}
\eta(dy\times\Box)&\defi& \int_{V}\int_{\RR^0_+} \delta_{\phi(x,u)}(dy) \II\{\phi(x,u)\in V\}\left( \int_{\bf A} p_T([u,\infty]|x,a)p_A(da|x)\right)~du~\mu(dx\times \bar\RR^0_+\times{\bf A})\nonumber\\
&=&\int_{\RR^0_+}\left\{\int_{V\times{\bf A}}\delta_{\phi(x,u)}(dy) \mu(dx\times[u,\infty]\times da)\right\} du\label{e10}\\
&=&\int_{\RR^0_+}\left\{\int_{\tilde V}\delta_{\phi((\tilde x,0),u)}(dy)\mu(d\tilde x\times\{0\}\times[u,\infty]\times{\bf A})\right\} du,\nonumber
\end{eqnarray}
we may write $
\int_{V\times \bar\RR^0_+\times \textbf{A}} \left\{ \int_{[0,\theta]}  C^g_j(\phi(x,u))  du
\right\}\mu(dx\times d\theta\times da)=\int_{V} C^g_j(y)\eta(dy\times\Box).
$

Similarly to the above, taking into account that the measure $\mu$ is concentrated on ${\bf M}\times{\bf A}$, we have that,
for each $j=0,1,\ldots,J$,
\begin{eqnarray*}
&&\int_{V\times \bar\RR^0_+\times \textbf{A}} \left\{ \II\{\theta<+\infty\}  C^I_j(\phi(x,\theta),a)\right\}
\mu(dx\times d\theta\times da)\\
&=&\int_{V\times \bar\RR^0_+\times \textbf{A}} \left\{ \II\{\theta<+\infty\}\II\{\phi(x,\theta)\in V\} C^I_j(\phi(x,\theta),a)\right\}
\mu(dx\times d\theta\times da)\\
&=& \int_{V\times{\bf A}} C^I_j(y,a) \eta(dy\times da),
\end{eqnarray*}
where
\begin{eqnarray}
\eta(dy\times da)&\defi&
 \int_{V}\int_{\RR^0_+} \delta_{\phi(x,\theta)}(dy) \II\{\phi(x,\theta)\in V\} \mu(dx\times d\theta\times da) \label{e11}
\end{eqnarray}
is a finite measure on $V\times{\bf A}$, since the measure $\mu$ is finite.

If Conditions  \ref{ConstrainedPPZcondition05},  \ref{ConstrainedPPZcondition01} and \ref{co31} are satisfied, and the space $\bf X$ is Borel, then the linear program (\ref{SashaLp02}) can now be rewritten as
\begin{eqnarray} \label{e106}
\mbox{Minimize}&:& \int_{V \times{\bf A}_\Box}{C}_0(y,a)\eta(dy\times da) \\
&&\mbox{over finite measures $\mu$ on $V\times{\bf B}=V\times\bar\RR^0_+\times{\bf A}$
concentrated on
${\bf M}\times{\bf A}$}\nonumber\\
\mbox{subject to}&:& \mbox{(\ref{e105}), (\ref{e10}), (\ref{e11}) and }
\int_{V \times{\bf A}_\Box}{C}_j(y,a)\eta(dy\times da)\le d_j,~j=1,2,\dots,J. \nonumber
\end{eqnarray}
The space ${\bf A}_\Box:={\bf A}\cup\{\Box\}$ and the functions $C_j$ are as introduced above: see (\ref{e103}). Proposition \ref{pr1} is valid for the linear program (\ref{e106}), too.

\begin{definition}\label{d101}
Suppose Conditions  \ref{ConstrainedPPZcondition05},  \ref{ConstrainedPPZcondition01} and \ref{co31} are satisfied, and assume that the space $\bf X$ is a Borel space.
For a finite measure $\mu$ on $V\times\bar\RR^0_+\times{\bf A}$ satisfying equation (\ref{e105}), the measure $\eta$ on $V\times{\bf A}_\Box$ defined by
\begin{equation}\label{e16p}
\eta(\Gamma_X\times\Gamma_A)\defi \eta(\Gamma_X\times (\Gamma_A\cap{\bf A}))+\eta(\Gamma_X\times\Box)\II\{\Box\in\Gamma_A\},~~~\Gamma_X\in{\cal B}(V),~\Gamma_A\in{\cal}({\bf A}_\Box),
\end{equation}
where the measures  $\eta(dy\times \Box)$ on $V$ and $\eta(dy\times da)$ on $V\times{\bf A}$ were introduced in (\ref{e10}) and (\ref{e11}),
is called the aggregated occupation measure (induced by $\mu$).
\end{definition}

In what follows, we will characterize  the aggregated measures $\eta$ without references to the measures $\mu$: see linear program (\ref{e123}).

\section{Main Results}\label{sec33}

\begin{definition}\label{d5}
We call the orbit of a point $\tilde x^0\in\tilde{\bf X}$ the following subset of ${\bf X}$:
$$_{\tilde x^0}{\cal X}=\{(\tilde\phi(\tilde x^0,t),t):~t\in\RR^0_+\}=\{\phi((\tilde x^0,0),t):~t\in\RR^0_+\}.$$
\end{definition}

We underline that the flow $\phi$ has no cycles and, if the flows $\phi$ and $\tilde\phi$ are continuous, then every orbit is a closed set in $\tilde{\bf X}\times\RR^0_+$.

\begin{condition}\label{co42}
Two different orbits  do not intersect, i.e., for any two distinct   points $\tilde{x}^0_1\ne\tilde{x}^0_2\in \tilde{{\bf X}}$,   $_{\tilde {x}^0_1}{\cal X}\cap{}_{\tilde {x}^0_2}{\cal X}=\emptyset$.
\end{condition}

\begin{definition}\label{d7}
Under Condition \ref{co42}, for each $y=(\tilde y,t)\in{\bf X}$, we introduce $h(y)$ equal to the point $\tilde x^0\in\tilde{\bf X}$ such that $\tilde y=\tilde\phi(\tilde x^0,t)$ and put $\tau_y = t$. The mappings $F:~\tilde{\bf X}\times \RR^0_+\to{\bf X}$ and $F^{-1}:~{\bf X}\to \tilde{\bf X}\times \RR^0_+$ are
defined as
\begin{eqnarray}\label{e102}
F(\tilde x^0,t)\defi (\tilde\phi(\tilde x^0,t),t)=\phi((\tilde x^0,0),t),~~~~~\mbox{ and }~F^{-1}(y)=(h(y),\tau_y).
\end{eqnarray}
\end{definition}

Note that the mapping $h:~{\bf X}\to \tilde{\bf X}$ is well defined: if, for $y=(\tilde y,t)\in{\bf X}$, for two points $\tilde x^0_1\ne\tilde x^0_2$ from $\tilde{\bf X}$, $\tilde y=\tilde\phi(\tilde x^0_1,t)=\tilde\phi(\tilde x^0_2,t)$, then the different orbits $_{\tilde x^0_1}{\cal X}$ and $_{\tilde x^0_2}{\cal X}$ are not disjoint having the common point $y$.

\begin{figure}[!htb]
	\centering
	\includegraphics[scale=0.3]{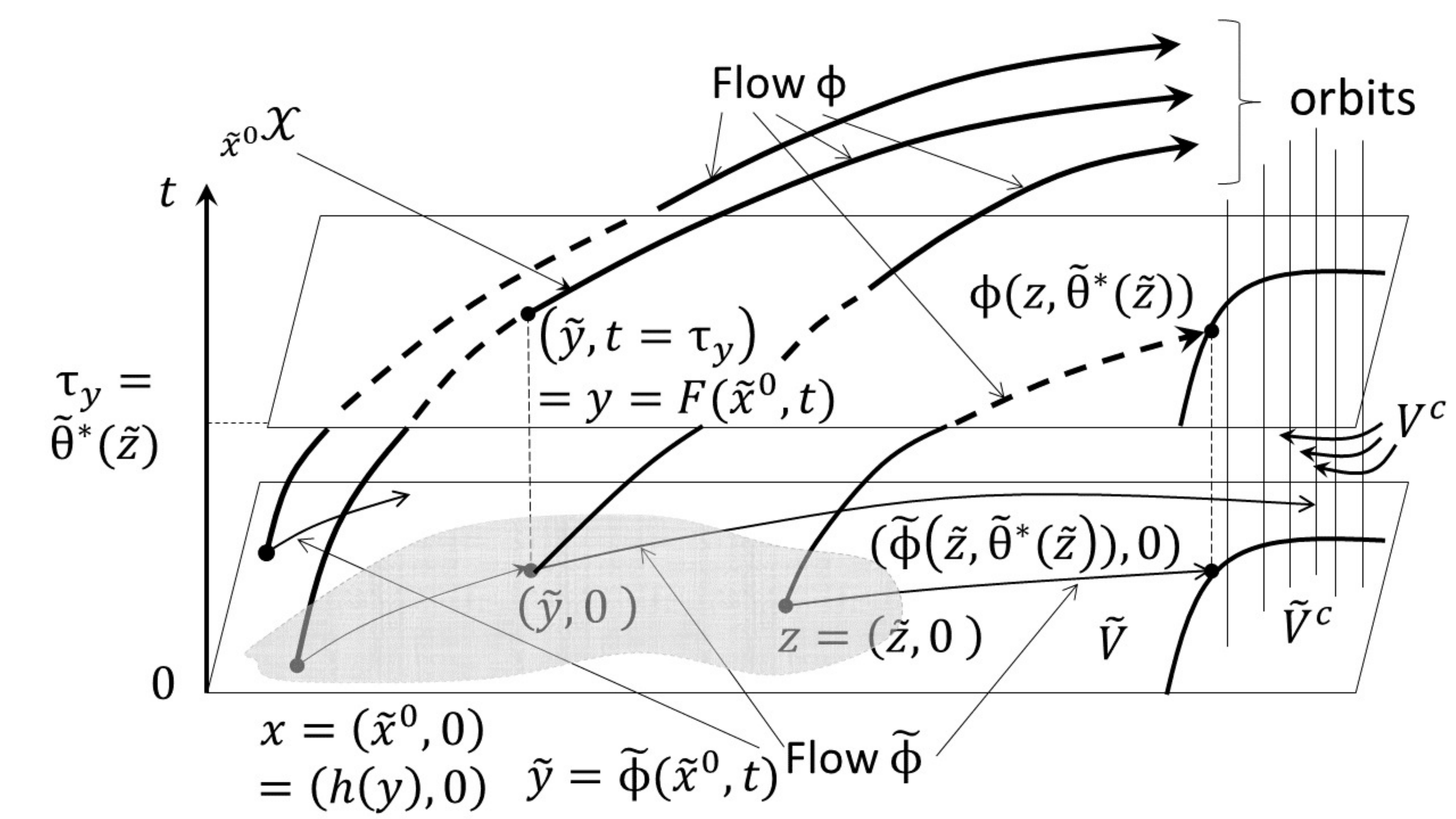}
	\caption{Flows $\tilde\phi$ and $\phi$. The grey area is $\tilde{\bf V}_\nu$: outside it $\hat\nu\equiv 0$.\protect\newline
	In general, the closed set $V^c$ can be arbitrary enough. Here, we assumed that the functions $C^g_j((\tilde x,t))$ do not depend on $t$, so that $V^c$ is the vertical cylinder.}
	\label{fig3}
\end{figure}

All the introduced notations are illustrated on Figure \ref{fig3}.

The mapping $F$ describes the forward movement from the starting point $(\tilde x^0,0)$ along the orbit $_{\tilde x^0}{\cal X}$; the inverse mapping $F^{-1}$ defines the starting point $\tilde x^0$, along with the duration of movement.

If Condition \ref{co42} is satisfied, one can define the flows $\tilde\phi$ and $\phi$ in the reverse time. For each $y=(\tilde y,t)\in{\bf X}$ we say that $\tilde\phi(\tilde y,-t)=h(y)$ and, for all $u\in[0,t]$, we put $\tilde\phi(\tilde y,-u)\defi \tilde\phi(h(y),t-u)$. For the flow $\phi$, we put $\phi(y,-u)=\phi((\tilde y,t),-u)=(\tilde\phi(\tilde y,-u),t-u)$.
The semigroup property here takes the form $\phi(x, t+s) = \phi(\phi(x,s), t)$ for $s$ and $t$ satisfying $s \ge -\tau_x,\ t+s \ge -\tau_x$. Note that $\tilde{\phi}$ in the reverse time is a function defined on $\{(y,-u):~y=(\tilde{y},t)\in \textbf{X},~0\le u\le t\}$.

The next condition requires that the speed of moving along the flow $\tilde\phi$ from $h((\tilde y,t))$ to $\tilde y$ is bounded.

\begin{condition} \label{con45}
Condition  \ref{co42} is satisfied, the  flows $\tilde{\phi}$ and $\phi$ are continuous, and there exists a $(0,\infty)$-valued function $d$ on $\mathbb{R}_+^0$, bounded on every finite interval $[0,T]$ and
such that for all $y_1=(\tilde y_1,t_1),~y_2=(\tilde y_2,t_2)\in{\bf X}$, \begin{eqnarray*}
\tilde\rho(h(y_1),h(y_2))=\tilde\rho(\tilde\phi(\tilde y_1,-t_1),\tilde\phi(\tilde y_2,-t_2))\le (d(t_1)\vee d(t_2))\rho(y_1,y_2),
\end{eqnarray*}
where $\rho$ and $\tilde{\rho}$ denote the compatible metrics on $\tilde{\bf X}\times\RR^0_+$ and $\tilde{\bf X}$, respectively.
\end{condition}

\begin{lemma}\label{la5}
Suppose Condition \ref{con45} is satisfied.
Then the mapping $h:~{\bf X}\to\tilde{\bf X}$, introduced in Definition \ref{d7},
is continuous, the flows $\tilde\phi$ and $\phi$ in the reverse time are continuous,  the mapping $F$
is a homeomorphism between $\tilde{\bf X}\times \mathbb{R}_+^0$ and ${\bf X}$, and the set $\bf X$ is a Borel space.
\end{lemma}

The proofs of this and several other auxiliary lemmas are postponed to the Appendix. Below, we assume that Condition \ref{con45} is satisfied.

For the points $\tilde x\in \tilde V$, the function $\tilde{\theta}^\ast(\tilde{x})$ defined by  (\ref{e14pp})
describes the time duration of the orbit $_{\tilde x}{\cal X}$ to be within the set $V$.
Recall that every orbit remains in $V^c$ after it reaches that set.

\begin{remark}\label{remark01}
Suppose \ref{ConstrainedPPZcondition05}, \ref{ConstrainedPPZcondition01}, \ref{co31}, and \ref{con45} are  satisfied. Then the mapping $F$ defined in (\ref{e102}) (its restriction on ${\bf D}$, to say more precisely) provides a homeomorphism and thus also an isomorphism between the sets
\begin{equation}\label{e15}
{\bf D}\defi\{(\tilde x^0,t):~\tilde x^0\in\tilde V,~0\le t<\tilde\theta^*(\tilde x^0)\}=\{(\tilde x^0,t):~\phi((\tilde x^0,0),t)\in  V\}
\end{equation}
and $V$.
Indeed, $F(\tilde x^0,t)\in V$ if and only if the pair $(\tilde x^0,t)$ belongs to the set $\bf D$. Thus, $F({\bf D})=V$ and $F^{-1}(V)={\bf D}$.
Recall that $\phi((\tilde x^0,0),t)\in V^c$ for all $t\ge\tilde\theta^*(\tilde x^0)$.
\end{remark}

We underline that the points $(\tilde x^0,t)\in{\bf D}$ and $(\tilde y,t)\in{\bf X}$ have different meanings, although the components $\tilde x^0,\tilde y\in\tilde{\bf X}$ and $t\in\RR^0_+$ look the same. That is the reason to equip the first coordinates of points in $\bf D$ with the upper index $0$, to make them look different from the points in $\bf X$. The pair $(\tilde x^0,t)\in{\bf D}$ is just the reference point of the orbit $_{\tilde x^0}{\cal X}$ and the duration of movement from $(\tilde x^0,0)$. It can easily happen that $(\tilde x^0,t)\notin{\bf X}$.

\begin{definition}\label{d8}
Suppose Conditions \ref{ConstrainedPPZcondition05}, \ref{ConstrainedPPZcondition01}, \ref{co31}, and \ref{con45} are satisfied. If $\zeta$ is a measure on $V$, then $\check{\zeta}$ denotes the image of $\zeta$ on $\bf D$ under the mapping $F^{-1}$:
\begin{eqnarray*}
\check{\zeta}(\Gamma)=\zeta(F(\Gamma)),~~~~~\Gamma\in{\cal B}({\bf D}).
\end{eqnarray*}
{In case the measure $\zeta$ is finite, we, with slight but convenient abuse of notations, introduce
$\hat\zeta(\tilde\Gamma)\defi\check{\zeta}(\tilde\Gamma\times \RR^0_+)$ for $\tilde\Gamma\in{\cal B}(\tilde V)$ and $\check{\zeta}(dt|\tilde x^0)$, the stochastic kernel} from $\tilde V$ to $\RR^0_+$
such that
\begin{eqnarray*}\check{\zeta}(d\tilde x^0\times dt)=\hat\zeta(d\tilde x^0)\check{\zeta}(dt|\tilde x^0),
\end{eqnarray*}
see \cite[Cor.7.27.2]{Bertsekas:1978} or \cite[Prop.D.8]{las96}.
\end{definition}

Clearly,
\begin{eqnarray*}
\zeta_1\le \zeta_2 ~\mbox{ set-wise } \Longleftrightarrow~ \check{\zeta}_1 \le \check{\zeta}_2 ~\mbox{ set-wise}.
\end{eqnarray*}

Note that if $\zeta$ is a finite measure, then $\check\zeta([0,\tilde\theta^*(\tilde x^0))|\tilde x^0)=1$ for $\hat\zeta$-almost all $\tilde x^0\in \tilde V$, and we extend the kernel $\check\zeta$ to $\RR^0_+$ by putting $\check\zeta([\tilde\theta^*(\tilde x^0),\infty)|\tilde x^0):=0$.
If the measure $\zeta$ is zero outside the set $\tilde V\times\{0\}$, then $\check{\zeta}(\Gamma)=0$ for all measurable subsets $\Gamma\subseteq {\bf D}\cap\{(\tilde x^0,t)\in\tilde{\bf X}\times \mathbb{R}_+^0:~t>0\}$, $\hat\zeta(\tilde\Gamma)=\zeta(\tilde\Gamma\times\{0\})$ for all $\tilde\Gamma\in{\cal B}(\tilde V)$, and $\check\zeta(dt|\tilde x^0)=\delta_0(dt)$ for $\hat\zeta$-almost all $\tilde x^0\in\tilde V$.

\begin{definition}\label{d9} Under Conditions \ref{ConstrainedPPZcondition05}, \ref{ConstrainedPPZcondition01}, \ref{co31} and \ref{con45},
a measure $\zeta$ on $V$ is called normal if there exist a finite measure $L$ on $\tilde V$ and a bounded measurable function $g(\tilde x^0,u):~\tilde V\times \RR^0_+\to\RR^0_+$ such that
$$\check\zeta (d\tilde x^0\times du)=g(\tilde x^0,u)du~L(d\tilde x^0).$$
Equivalently, for all $\Gamma_X\times\Gamma_t\in {\cal B}(V)$,
\begin{eqnarray*}
\zeta(\Gamma_X\times\Gamma_t)&=&\int_{\bf D}\II\{F(\tilde x^0,u)\in\Gamma_X\times\Gamma_t\}\check\zeta(d\tilde x^0\times du)\\
&=&\int_{\tilde V}\int_{[0,\tilde\theta^*(\tilde x^0))}\delta_{\tilde\phi(\tilde x^0,u)}(\Gamma_X)\delta_u(\Gamma_t) g(\tilde x^0,u)du~L(d\tilde x^0).
\end{eqnarray*}
See Remark \ref{remark01}.

A measure $\eta$ on $V\times{\bf A}_\Box$ is called normal if $\eta(V\times{\bf A})<\infty$ and the measure $\eta(dx\times\Box)$ on $V$ is normal.
\end{definition}

Clearly, every normal measure is $\sigma$-finite. Similarly, a normal measure $\zeta$ defined on some orbit  $_{\tilde z}{\cal X}\subseteq \textbf{X}$ is understood:
\begin{eqnarray*}
\check\zeta(d\tilde x^0\times du)=g(u)du~\delta_{\tilde z}(d\tilde x^0).
\end{eqnarray*}

\begin{lemma}\label{l101}
Suppose Conditions \ref{ConstrainedPPZcondition05}, \ref{ConstrainedPPZcondition01}, \ref{co31} and \ref{con45} are satisfied. Then the following assertions hold true.
\begin{itemize}
\item[(a)] For every finite measure $\mu$ on $V\times\bar\RR^0_+\times{\bf A}$ satisfying equality (\ref{e105}), the induced aggregated occupation measure $\eta$ on $V\times{\bf A}_\Box$ is normal.
\item[(b)] If $\eta^1$ and $\eta^2$ are two normal measures on $V\times{\bf A}_\Box$ such that $\eta^1\ge \eta^2$ set-wise, and thus the difference $\eta^1-\eta^2$ is a (positive) measure, then $\eta:=\eta^1-\eta^2$ is also a normal measure on $V\times{\bf A}_\Box$.
\end{itemize}
\end{lemma}

In Definition \ref{d3}, we introduce the class of so called test functions used to characterize measures on $\bf X$.

\begin{definition}\label{d3}
$\bf W$ is the space of measurable bounded functions $w$ on $\bf X$, absolutely continuous, either negative and increasing or positive and decreasing along the flow  $\phi$ (see Definition \ref{JulyDef01}) and satisfying conditions
\begin{itemize}
\item $w(y)=0$ for all $y\in V^c$ and
\item $\lim_{t\to\infty} w(\phi(x,t))=0$ for all $x\in V$ such that  $\phi(x,t)\in V$ for all $t\in\mathbb{R}_+^0$.
\end{itemize}
\end{definition}

Throughout this paper, $\chi w$ denotes a function as in Lemma \ref{l1} (see Appendix).
Without loss of generality, one can assume, for each negative (or positive) function $w\in \textbf{W}$, that the function $\chi w$ is  {positive} (or  {negative}), i.e., in (\ref{ea1}) one can put $g(\cdot)\equiv 0$. Note that below we consider only such measures $\zeta$ on $V$ that the value of the integral $\int_{V}\chi w(x)\zeta(dx)$ does not depend on the function $g$ in (\ref{ea1}).

Suppose Conditions \ref{ConstrainedPPZcondition05}, \ref{ConstrainedPPZcondition01}, \ref{co31}, and \ref{con45} are satisfied and introduce the following linear program
\begin{eqnarray}
\mbox{Minimize over}\nonumber \\
\mbox{the normal measures $\eta$ on $V\times \textbf{A}_\Box$}&:&\int_{V \times\textbf{A}_\Box}{C}_0(x,a)\eta(dx\times da) \label{e123} \\
\mbox{subject to}&:&
w(x_0)+\int_{V} \chi w(x)\eta(dx\times\Box)-\int_{V} w(x)\eta(dx\times{\bf A})\label{e17}\\
&&+\int_{V\times{\bf A}} w(l(x,a))\eta(dx\times da)=0~~~~~\forall~w\in{\bf W};\nonumber \\
&&\int_{V \times\textbf{A}_\Box}{C}_j(x,a)\eta(dx\times da)\le d_j,~j=1,2,\dots,J. \nonumber
\end{eqnarray}
The space ${\bf A}_\Box$ and functions $C_j$ were defined in Section \ref{sec3}: see (\ref{e103}).

\begin{remark}\label{rem991}
Compared with (\ref{SashaLp02}) and (\ref{e106}), the dimensionality of the linear program (\ref{e123}),(\ref{e17}) is reduced in the sense that the measures $\mu$ were on the space $V\times\bar\RR^0_+\times{\bf A}$, and the measures $\eta$ are on the space $V\times({\bf A}\cup\{\Box\})$. Therefore, e.g., from the computational point of view, the linear program (\ref{e123}),(\ref{e17}) is easier.
\end{remark}

We are ready to formulate the main results.

\begin{theorem}\label{t1}
Suppose Conditions \ref{ConstrainedPPZcondition05}, \ref{ConstrainedPPZcondition01}, \ref{co31} and \ref{con45} are satisfied. Then, for every finite measure $\mu$ on $V\times\bar\RR^0_+\times{\bf A}$, concentrated on ${\bf M}\times{\bf A}$ and satisfying equality (\ref{e105}), its induced aggregated occupation measure $\eta$ on $V\times{\bf A}_\Box$
satisfies equation (\ref{e17}) for all functions $w\in{\bf W}$. All the integrals in (\ref{e17}) are finite.
\end{theorem}
The proof of this statement is postponed to Section \ref{sec5}.

\begin{theorem}\label{t2} Suppose Conditions \ref{ConstrainedPPZcondition05}, \ref{ConstrainedPPZcondition01}, \ref{co31} and \ref{con45} are satisfied. Then
every normal measure $\eta$ on $V\times{\bf A}_\Box$, satisfying equation (\ref{e17}),
uniquely defines a reasonable Markov strategy $\pi^\eta$ (called ``induced'' by $\eta$)
such that, for the aggregated occupation measure $\tilde\eta$ defined by (\ref{e16p}) (recall  (\ref{e10}) and (\ref{e11})) with $\mu$ being replaced by the occupation measure $\mu^{\pi^\eta}$ of the strategy $\pi^\eta$ as in (\ref{e101}) with $\pi=\pi^\eta$, the following inequalities hold:
\begin{eqnarray*}
\tilde\eta(\Gamma)\le\eta(\Gamma)~\forall~\Gamma\in{\cal B}(V\times{\bf A}_\Box).
\end{eqnarray*}
\end{theorem}
The proofs of Theorem \ref{t2}  and of the next corollary are postponed to Section \ref{sec6}.

\begin{corollary}\label{corol1}
Let Conditions \ref{ConstrainedPPZcondition05}, \ref{ConstrainedPPZcondition01}, \ref{co31}, and \ref{con45} be satisfied. Then linear program (\ref{e106}) is equivalent to linear  program (\ref{e123}).

To be more precise, if the finite measure $\mu^*$ on $V\times\bar\RR^0_+\times{\bf A}$ solves linear program (\ref{e106}), then the measure $\eta^*$ on $V\times{\bf A}_\Box$, given by (\ref{e10}), (\ref{e11}) and (\ref{e16p}), i.e., the aggregated occupation measure induced by $\mu^*$, solves linear program (\ref{e123}). Conversely, if the measure $\eta^*$ on $V\times{\bf A}_\Box$ solves linear program (\ref{e123}), then, for the Markov strategy $\pi^*$ induced by $\eta^*$ as in Theorem \ref{t2},  the corresponding occupation measure $\mu^{\pi^*}$ on $V\times\bar\RR^0_+\times{\bf A}$, defined in (\ref{e101}), solves linear program (\ref{e106}).
\end{corollary}

According to Corollary \ref{corol1} and Section \ref{sec3} (see Proposition \ref{pr1}), the minimal values of the linear programs (\ref{SashaLp02}) and (\ref{e123}) coincide and equal  the minimal value of the original problem (\ref{PZZeqn02}).
As soon as the optimal solution $\eta^*$ to the  linear program (\ref{e123})  is obtained, the induced Markov strategy $\pi^*$,  solves the original optimal impulsive control problem (\ref{PZZeqn02}): see Proposition \ref{pr1} and remember that the linear programs (\ref{SashaLp02}) and (\ref{e106}) are  equivalent.
Recall that linear program (\ref{e106}) has an optimal solution by Proposition \ref{pr1}; hence the linear program (\ref{e123}) is also solvable. Note also that, having in hand the Markov strategy $\pi^*$, one can compute the corresponding occupation measure $\mu^{\pi^*}$ (\ref{e101}), and after that the stationary strategy as in Proposition \ref{pr1}  also solves the optimal impulsive control problem (\ref{PZZeqn02}).

For the discussions in the rest of this section, we suppose all the mappings and functions $l,C^g_j$, and $C^I_j$ do not depend on the component $t$ of the state $x=(\tilde x,t)$. Then the linear program (\ref{e123}) is actually in terms of (marginal) measures $\tilde\eta(d\tilde x\times da)$ on $\tilde V\times{\bf A}_\Box$ defined by
$$\tilde\eta(\Gamma_X\times\Gamma_A)\defi\eta([(\Gamma_X\times[0,\infty))\cap V]\times\Gamma_A)=\int_{(\Gamma_X\times[0,\infty))\cap V} \eta(d\tilde x\times dt\times\Gamma_A).$$
The marginals $\tilde\eta$ of normal measures $\eta$ (naturally called normal on $\tilde V\times{\bf A}_\Box$) are characterized  as follows: $\tilde\eta(\tilde V\times{\bf A})<\infty$ and there exist a finite measure  $L$ on $\tilde V$ and a bounded non-negative measurable function $g$ on $\tilde V\times\RR^0_+$ such that
$$\tilde\eta(\Gamma_X\times\Box)=\eta([(\Gamma_X\times[0,\infty))\cap V]\times\Box)=\int_{\tilde V}\int_{[0,\tilde\theta^*(\tilde x))} \delta_{\tilde\phi(\tilde x,u)}(\Gamma_X) g(\tilde x,u)du~L(d\tilde x)$$
(see Definition \ref{d9}). The  test functions $\tilde w$ on $\tilde{\bf X}$ are
measurable bounded, absolutely  continuous, either negative and increasing or positive and decreasing along the flow, and
such that $\tilde w(y)=0$ for all $y\in\tilde V^c$ and $\lim_{t\to\infty}\tilde w(\tilde\phi(\tilde x,t))=0$ for all $\tilde x\in\tilde V$ such that $\tilde\phi(\tilde x,t)\in\tilde V$ for all $t\in\RR^0_+$.

This linear program, in terms of marginal measures $\tilde\eta$, is solvable under Conditions \ref{ConstrainedPPZcondition05}, \ref{ConstrainedPPZcondition01}, \ref{co31}, and \ref{con45}. (The last condition is
for the model with the extended state space ${\bf X}\subset\tilde{\bf X}\times\RR^0_+$.) The minimal value of this program coincides with the minimal value of the original  problem (\ref{PZZeqn02}). Therefore, when reformulating the  optimal impulsive control problem in terms of aggregated occupation measures, the extension of the state space, as in Remark \ref{remark2}, is not needed. 

On the other hand, the construction of the optimal control strategy $\pi^*$, induced by the optimal solution $\eta^*$ to the linear program (\ref{e123}), is essentially based on the analysis in the extended state space: see the proof of Theorem \ref{t2}. Note that, in the case of the extended state space, the (full) orbits in $\bf X$, as in Definition \ref{d5}, form a Borel space because they are characterized by the starting points $\tilde x^0$. If one manages to describe the space of orbits in $\tilde{\bf X}$ as a Borel space, then one can avoid such an extension of the basic state space $\tilde{\bf X}$.

\section{Example and Comparison with Other Works} \label{sec7}

Consider the following simple but not trivial optimal impulse control problem in the space $\RR^0_+$.
\begin{equation}\label{em2}
\left.\begin{array}{rcl}
d\tilde x &:=& G(\tilde x) dt+dW(t),~~~~~\tilde x(0-)=\tilde x_0>0;\\~&~&~\\
&& \displaystyle \int_0^\infty C^g(\tilde x(u))du+\int_0^\infty dW(u) \to \inf_W,
\end{array}\right\}
\end{equation}
where
$$W(u)\defi \sum_{j=1}^\infty\left(\sum_{i=1}^{j-1} a_i\right) I\{T_{j-1}\le u<T_j\};$$
$$0=T_0\le T_1\le T_2\le \ldots,~~~T_j\in[0,\infty];~~T_{j-1}=T_j~\mbox{only if } T_{j-1}=\infty;~~~\lim_{j\to\infty}T_j=\infty.$$
The impulse control strategy $W$, represented by $\{T_j,a_j\}_{j=1}^\infty$, can be arbitrary, satisfying the condition $a_j\ge\delta>0$.
The measurable functions $G(\cdot)>\delta\ge 0$ and $C^g(\cdot)\ge 0$ are fixed and smooth enough, such that $H\defi\int_0^\infty C^g(\tilde X(u))du<\infty$. Here $\tilde X(\cdot)$ is the solution to (\ref{em2}) when $W(u)\equiv 0\Leftrightarrow T_1=\infty$.

Clearly, this problem can be easily reformulated in terms of Section \ref{sec2}. The flow $\tilde\phi$ on $\tilde{\bf X}=[\tilde x_0,\infty)$ comes from the differential equation (\ref{em2}) at $W(\cdot)\equiv 0$; ${\bf A}=[\delta,H]$ (no reason to apply the impulses $a>H$); $l(x,a)=\tilde l((\tilde x,t),a)=\tilde l(\tilde x,a)=\tilde x+a$; as usual, $t$ is the time elapsed since the most recent impulse. We consider the unconstrained case with $J=0$. The gradual cost rate is $C^g(\cdot)$, and the cost of the impulse $a\in{\bf A}$ equals $C^I(x,a)=a$. Simultaneous impulses of the sizes $a,b,\ldots$ can be considered as one impulse of the size $a+b+\ldots$. We assume that, for some $K\in(\tilde x_0,\infty)$, $C^g(\tilde x)=0$ for $\tilde x\ge K$ and $C^g(\tilde x)>0$ for $\tilde x<K$, so that $\tilde V=[0,K)$, and $\tilde V^c=[K,\infty)$.
Now, all the Conditions \ref{ConstrainedPPZcondition05}, \ref{ConstrainedPPZcondition01}, \ref{co31}, and \ref{con45} are satisfied, and hence the minimal value of the impulse control problem (\ref{em2}) coincides with the minimal value of the linear program (\ref{e123}).

\begin{figure}[!htb]
	\centering
	\includegraphics[scale=0.3]{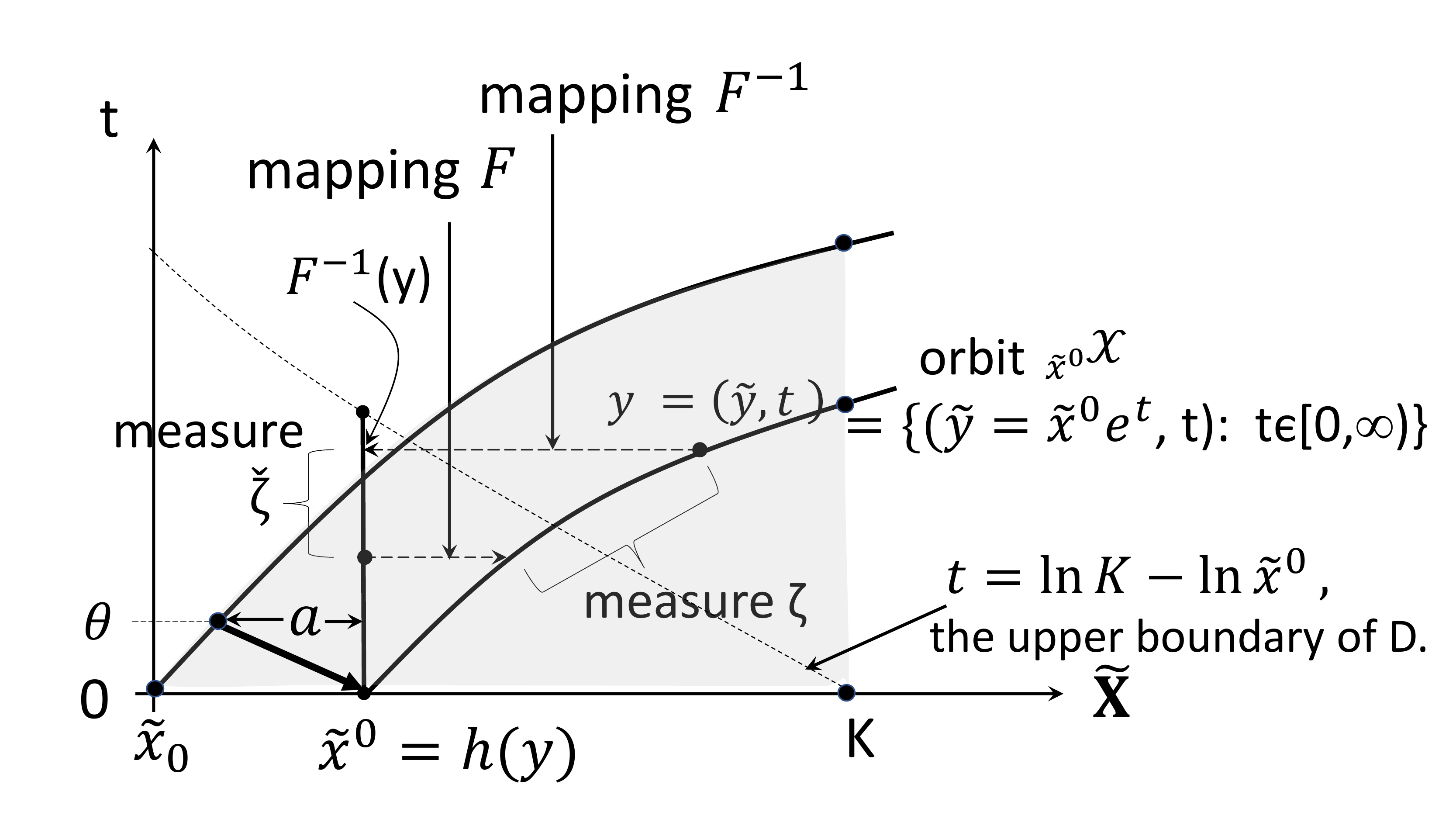}
	\caption{$\bf X$ is the area below the orbit starting from $\tilde x_0$;
	the grey area is $V$. The bold arrow leading to the point $\tilde x^0$ represents the impulse of the size $a$ applied at the time moment $T_1=\theta$.}
		\label{fig8}
\end{figure}

One also can illustrate all the definitions introduced in Section \ref{sec33}. To be specific, take $G(\tilde x)=\tilde x$, so that $\tilde\phi(\tilde x,t)=\tilde xe^t$ is the solution to (\ref{em2}) starting from $\tilde x$, when $W(u)\equiv 0$, i.e., $T_1=\infty$;
\begin{eqnarray*}
{\bf X}=\{(\tilde y,t):~\tilde y\in\tilde{\bf X}=[\tilde x_0,\infty),~ t\in[0,\ln\frac{\tilde y}{\tilde x_0}]\}:
\end{eqnarray*}
a point $\tilde y$ cannot appear later than $\ln\frac{\tilde y}{\tilde x_0}$ time units after any one impulse. See Figure \ref{fig8}. Now
$$h(y)=h(\tilde y,t)=\frac{\tilde y}{e^t};~F(\tilde x^0,t)=(\tilde x^0 e^t,t);~F^{-1}(y)=F^{-1}(\tilde y,t)=(\frac{\tilde y}{e^t},t);~~~\tilde\theta^*(\tilde x)=\ln \frac{K}{\tilde x}~\mbox{ for } \tilde x\in V;$$
$$V=\{(\tilde y,t):~\tilde y\in\tilde V=[0,K),~t\in[0,\ln\frac{\tilde y}{\tilde x_0}]\};~
{\bf D}=\{(\tilde x^0,t):~\tilde x^0\in[0,K),~0\le t<\ln\frac{K}{\tilde x^0}\}.$$
With some abuse of notations, we avoid the double brackets in the expressions like $h(y)=h((\tilde y,t))$.
The mappings $F$ and $F^{-1}$ are one-to-one and continuous. According to Definition \ref{d9}, a measure $\zeta$ on $V$ is normal if and only if the conditional distribution $\check\zeta(dt|\tilde x^0)$ is ($L$-almost surely) absolutely continuous with respect to the Lebesgue measure, that is, the measure $\zeta$, restricted to the orbit $~_{\tilde x^0}\cal X$, is ($L$-almost surely) absolutely continuous with respect to the Lebesgue measure on that orbit.

According to the last paragraphs in Section \ref{sec33}, we formulate the linear program (\ref{e123}) in terms of the (normal) marginal measures $\tilde\eta(d\tilde x\times da)$ on $\tilde V\times{\bf A}_\Box$, ignoring the $t$ component, time elapsed since the most recent impulse. The unnecessary `tilde' is omitted up to the end of this section, apart from $\tilde x_0$ (initial state).
\begin{equation}\label{e1235}
\left.\begin{array}{rcl}
\mbox{Minimize over the normal} \\
\mbox{measures $\eta$ on $[\tilde x_0,K)\times \textbf{A}_\Box$}&:&\displaystyle \int_{[\tilde x_0,K) }{C}^g( x)\eta(d x\times \Box)+\int_{[\tilde x_0,K)\times{\bf A}} a~\eta(d x\times da)  \\~\\
\mbox{subject to}&:&
w(\tilde x_0)+\int_{[\tilde x_0,K)} \chi w(x)\eta(dx\times\Box)\\
&&\displaystyle+\int_{[\tilde x_0,K)\times{\bf A}} [w(( x+a)-w(x)]\eta(d x\times da)=0.\end{array}\right\}
\end{equation}
The test functions $w$ on $[\tilde x_0,\infty)$ are bounded, measurable, absolutely continuous, either negative and increasing or positive and decreasing, and such that $w(x)=0$ on $[K,\infty)$. The measures $\eta(dx\times\Box)$ are finite on $[\tilde x_0,K)$ because $\theta^*(x)\le K/\delta$: recall that $G(\cdot)>\delta$; $\eta([\tilde x_0,K)\times{\bf A})<\infty$ due to the definition of a normal measure.

It is interesting to compare the linear program  (\ref{e1235}) with the linear programs which appeared in \cite{b13,b13p,b13pp}. In those articles, the impulse control problem was formulated on the finite time horizon $[0,T]$, but the constructions can be formally adjusted for $[0,\infty)$.

Following the ideas of \cite{b13}, the problem (\ref{em2}) is replaced with  the following linear program on the space of the so called occupation measures $\Upsilon^1_1$ and $\Upsilon^1_2$:
\begin{equation}\label{e1236}\left.\begin{array}{rcl}
\mbox{Minimize over the finite}
\\ \mbox{measures $\Upsilon^1_1$ and $\Upsilon^1_2$ on $[\tilde x_0,\infty)$}
&:&\displaystyle\int_{[\tilde x_0,\infty) }{C}^g( x)\Upsilon^1_1(d x)+\int_{[\tilde x_0,\infty)}\Upsilon^1_2(d x)  \\~\\
\mbox{subject to}&:& \displaystyle
w(\tilde x_0)+\int_{[\tilde x_0,\infty)} \frac{dw(x)}{dx}G(x)\Upsilon^1_1(dx)\\~\\
&&\displaystyle+\int_{[\tilde x_0,\infty)} \frac{dw(x)}{dx}\Upsilon^1_2(d x)=0,
\end{array}\right\}
\end{equation}
where the test functions $w$ are continuously differentiable on $[\tilde x_0,\infty)$ and $\lim_{x\to\infty} w(x)=0$.

Consider the test functions $w$ as in (\ref{e1235}), which are continuously differentiable on $[\tilde x_0,\infty)$. Now $\chi w(x)=\frac{dw}{dx}G(x)$ and, for the measures
\begin{eqnarray*}
&&\Upsilon^1_1(dx):=\left\{\begin{array}{ll}
\eta(dx\times\Box), & \mbox{if } dx\in{\cal B}([\tilde x_0,K));\\
\mbox{arbitrary} & \mbox{otherwise}, \end{array}\right.\\
&\mbox{and}& \Upsilon^1_2(dz):=\left[\int\limits_{[\tilde x_0,K)\times{\bf A}} \hspace{-2mm}\II\{x<z<x+a\}\eta(dx\times da)\right] dz,
\end{eqnarray*}
all the expressions in (\ref{e1236}), take the form of those in (\ref{e1235}) because
\begin{eqnarray*}
\int_{[\tilde x_0,\infty)}\Upsilon^1_2(dz)&=& \int_{[\tilde x_0,K)\times{\bf A}}\int_{[\tilde x_0,\infty)}\II\{x<z<x+a\}~dz~\eta(dx\times da)\\&=&
\int_{[\tilde x_0,K)\times{\bf A}}a~\eta(dx\times da);\\
\int_{[\tilde x_0,\infty)}\frac{dw(z)}{dz}\Upsilon^1_2(dz)&=& \int_{[\tilde x_0,K)\times{\bf A}}\int_{[\tilde x_0,\infty)}\frac{dw(z)}{dz}\II\{x<z<x+a\}~dz~\eta(dx\times da)\\&=&
\int_{[\tilde x_0,K)\times{\bf A}}[w(x+a)-w(x)]\eta(dx\times da).
\end{eqnarray*}
Recall also that $C^g(x)=\frac{dw(x)}{dx}=0$ for $x\in[K,\infty)$.

In the works \cite{b13p,b13pp}, the impulse control problem (\ref{em2}) is formulated in a different way which is briefly presented below. The generic notations of \cite{b13p,b13pp} are changed to avoid the confusion with the notations in the present paper. Let a reasonable deterministic stationary control strategy, defined by $\{T_j,a_j\}_{j=1}^\infty$ and denoted below as $f:~{\bf X}_\Delta\to\bar\RR^0_+\times{\bf A}$, be fixed, such that $T_j=\infty$ if $x(T_j-)\ge K$. By the way, the number of finite moments $T_j$ is finite, and the class of such strategies is sufficient in the unconstrained problem (\ref{em2}) by Theorem 1 in \cite{PiunovskiySasha:2018}. Introduce the measure
$$\sigma(dt):= dt+\sum_{j=1}^\infty \II\{T_j<\infty\}\delta_{T_j}(dt)$$
on the time scale $[0,\infty)$. The model (\ref{em2}) is represented as
\begin{equation}\label{est1}
\left.\begin{array}{rcl}
dx(t)&=&\displaystyle\int_{{\bf A}\cup\{0\}}\frac{ F(x(t),a)}{1+|a|^p}\kappa(da|t)\sigma(dt);~~x(0-)=\tilde x_0;\\~&~&~\\
&&\displaystyle\int_{[0,\infty)}\int_{{\bf A}\cup\{0\}}\frac{L(x(t),a)}{1+|a|^p}\kappa(da|t)\sigma(dt)\to \inf_{\sigma,\kappa}\end{array}\right\}
\end{equation}
with the following system primitives:
\begin{itemize}
\item $p\ge 1$ is some fixed natural number.
\item $\kappa(da|t)=\delta_0(da)$ if $t$ is different from all $T_j$, so that $F(x,0)=G(x)$, and $a=0$ corresponds to the  absence of impulses.
\item $\kappa(da|T_j)=\delta_{a_j}(da)$ for $T_j<\infty$, and at that time moment the following fictitious process is introduced:
$$\frac{dy_{T_j}^{a_j}(u)}{du} =\frac{F(y_{T_j}^{a_j}(u),a_j)}{1+|a_j|^p},~~~y_{T_j}^{a_j}(0)=x(T_j-),~~~u\in[0,1].$$ The form of the function $F$ is seen in the next item.
\item  $x(T_j+)=x(T_j-)+[y^{a_j}_{T_j}(1)-y^{a_j}_{T_j}(0)]=
y_{T_j}^{a_j}(1)$.  To be consistent with the model (\ref{em2}), we should have $x(T_j+)=x(T_j-)+a_j$, so that for $a\in{\bf A}$ we put $F(y,a):=a(1+|a|^p)$.
\item Similarly, for consistency, we put $L(x,0)=C^g(x)$ and $L(y,a):=a(1+|a|^p)$ for $a\in{\bf A}$.
\end{itemize}
The occupation measure  on $[0,\infty)\times[\tilde x_0,\infty)\times[{\bf A}\cup\{0\}]$ as in \cite{b13p,b13pp},  corresponding to the strategy $f$ (equivalently, to the pair  $(\sigma,\kappa$)),  equals
$$\Upsilon^f(dt\times dy\times da):=\Phi(dy|a,t)\kappa(da|t)\sigma(dt)= \Upsilon^f_1(dt\times dy)\delta_0(da)+\Upsilon^f_2(dt\times dy\times da),$$
where
$$\Phi(dy|a,t):=\left\{\begin{array}{ll}
\delta_{x(t)}(dy), & \mbox{ if } t\ne T_j \mbox{ for all } j=1,2,\ldots;\\~\\
\int_0^1\delta_{y^a_t(u)}(dy)du=\int_0^1 \delta_{x(T_j-)+au}(dy) du, & \mbox{ if } t=T_j;
\end{array}\right.$$
$x(\cdot)$ is the trajectory of the system (\ref{est1}) (equivalently, of the system (\ref{em2})) under the strategy $f$. The presentation $\Upsilon^f=\Upsilon^f_1+\Upsilon^f_2$ corresponds to the decomposition of the measure $\sigma$ to the absolutely continuous and discrete parts.
Different control strategies $f$ as above, that is, different pairs $(\sigma,\kappa)$ define all different measures $\Upsilon$ under consideration,  which are denoted below as $\Upsilon^2$.

Below, the test functions $w$ are as in (\ref{e1235}) and continuously differentiable on $[\tilde x_0,\infty)$. In the linear program for the problem (\ref{est1}), suggested in \cite{b13p,b13pp}, all the integrated functions do not depend on time $t$. Thus, we immediately introduce the marginals $\hat\Upsilon^2(dy\times da):=\int_{[0,\infty)}\Upsilon^2(dt\times dy\times da)$:
\begin{eqnarray*}
\hat\Upsilon^2_1(dy)&:=& \int_{[0,\infty)}\Upsilon^2_1(dt\times dy)=\int_{[0,\infty)}\delta_{x(t)}(dy)~dt;\\
\hat\Upsilon^2_2(dy\times da)&:=& \int_{[0,\infty)}\Upsilon^2_2(dt\times dy\times da)=\sum_{T_j}\II\{T_j<\infty\}\left[\int_0^1\delta_{x(T_j-)+au}(dy)~du\right]\delta_{a_j}(da).
\end{eqnarray*}
Here the measure $\Upsilon^2=\Upsilon^2_1+\Upsilon^2_2$ comes from the pair $(\sigma,\kappa)$, which also defines the trajectory $x(\cdot)$ of the system (\ref{est1}); the measure $\hat\Upsilon^2_2$ is finite and $\hat\Upsilon^2_2([K+H,\infty)\times {\bf A})=0$ because $x(T_j-)<K$ and $a_j\le H$. The linear program as in \cite{b13p,b13pp} has the form
$$\left.\begin{array}{rcl}
\mbox{Minimize over}\\
\mbox{the measures $\hat\Upsilon^2$ on $[\tilde x_0,\infty)\times{\bf A}$}&:& \displaystyle \int_{[\tilde x_0,\infty)\times{\bf A}}\frac{L(y,a)}{1+|a|^p}~\hat\Upsilon^2(d y\times da) \\~\\
\mbox{subject to}&:& \displaystyle w(\tilde x_0)+ \int_{[\tilde x_0,\infty)\times{\bf A}}\frac{dw(y)}{dy}\frac{F(y,a)}{1+|a|^p}~\hat\Upsilon^2(d y\times da)=0 \end{array}\right\}$$
(note that $\lim_{y\to\infty}w(y)=0$),
or, more explicitly,
\begin{equation}\label{e1237}\left.\begin{array}{rcl}
\mbox{Minimize over}
\\ \mbox{the finite measures $\hat\Upsilon^2_1$ and $\hat\Upsilon^2_2$}
&:&\displaystyle\int_{[\tilde x_0,K) }{C}^g( y)\hat\Upsilon^2_1(d y)+\int_{[\tilde x_0,\infty)\times{\bf A}}a~\hat\Upsilon^2_2(d y\times da)  \\~\\
\mbox{subject to}&:& \displaystyle
w(\tilde x_0)+\int_{[\tilde x_0,K)} \frac{dw(y)}{dy}G(y)\hat\Upsilon^2_1(dy)\\~\\
&&\displaystyle+\int_{[\tilde x_0,\infty)\times{\bf A}} \frac{dw(y)}{dy}~a~\hat\Upsilon^2_2(d y\times da)=0.
\end{array}\right\}
\end{equation}
We underline that the measure $\hat\Upsilon^2_1$ is of no importance on $[K,\infty)$ because there $C^g(y)=\frac{dw(y)}{dy}=0$; it is finite on $[\tilde x_0,K)$ because $\int_{[0,\infty)}\delta_{x(t)}([\tilde x_0,K))dt\le\theta^*(\tilde x_0)<\infty$.

The measures $\hat\Upsilon^2_1$ and $\hat\Upsilon^2_2$ can be calculated based on the measures $\eta$ in (\ref{e1235}), so that all the expressions in (\ref{e1237}) become equal to those in (\ref{e1235}). Indeed, we put $\hat\Upsilon^2_1:=\Upsilon^1_1$ and
$$\hat\Upsilon^2_2(dy\times da):=\frac{1}{a}\left[\int_{[\tilde x_0,K)}\II\{x<y<x+a\}\eta(dx\times da)\right]~dy.$$
Now $\displaystyle \int_{\bf A} a~\hat\Upsilon^2_2(dy\times da)=\Upsilon^1_2(dy)$ and all the expressions in (\ref{e1237})  coincide with those in (\ref{e1236}) and, as shown above, are equal to those in (\ref{e1235}).

\section{Proof of Theorem \ref{t1}}\label{sec5}

\underline{Proof of Theorem \ref{t1}}. Note that, for each function $w\in{\bf W}$, for each fixed $x\in V$, the function $w(\phi(x,\cdot))$ is bounded on $\RR^0_+$.

According to Lemma \ref{l1}, for each fixed $x\in V$,
$$w(\phi(x,\theta))=w(x)+\int_{[0,\theta]}\chi w(\phi(x,s))ds,$$
where the function $\chi w$ is given by (\ref{ea1}).
After we integrate this equation over $V\times\RR^0_+$ with respect to the measure
$$\int_{\bf A} \II\{\phi(x,\theta)\in V\} p_T(d\theta|x,a)p_A(da|x)\mu(dx\times\bar\RR^0_+\times{\bf A}),$$
on $V\times\RR^0_+$,
where  the stochastic kernels $p_T$ and $p_A$ are as in (\ref{e14p}), we obtain the equality
\begin{eqnarray*}
\int_{V} w(y) \eta(dy\times{\bf A}) &=& \int_{V}\int_{\RR^0_+} w(\phi(x,\theta)) \II\{\phi(x,\theta)\in V\} \hat p(d\theta|x)\mu(dx\times\bar\RR^0_+\times{\bf A})\\
&=& \int_{V} w(x)\int_{\RR^0_+} \II\{\phi(x,\theta)\in V\} \hat p(d\theta|x) \mu(dx\times\bar\RR^0_+\times{\bf A})\\
&&+\int_{V}\int_{\RR^0_+}\II\{\phi(x,\theta)\in V\}\int_{[0,\theta]} \chi w(\phi(x,s)) ds~\hat p(d\theta|x)
 \mu(dx\times\bar\RR^0_+\times{\bf A}),
\end{eqnarray*}
where $\hat p(d\theta|x)\defi \int_{\bf A} p_T(d\theta|x,a)p_A(da|x)$. Note that all the integrals here are finite because the function $w(\cdot)$ is bounded and the measures $\mu$ and $\eta(dy\times{\bf A})$ are finite. For each $x\in V$, let us denote
$$\theta^*(x)\defi \inf\{\theta\in\RR^0_+:~\phi(x,\theta)\in V^c\}.$$
As usual, $\inf \emptyset \defi +\infty$. Since the flow $\phi$ is continuous, the function $\theta^*(\cdot)$ is measurable: see \cite[Lemma 27.1]{b1} or \cite[Prop.1.5, p.154]{b3}. Besides, $\theta^*(x)>0$ because the set $V$ is open and the set $V^c$ is closed.

Since the set $V^c$ is closed and the flow $\phi$ is continuous, in case $\theta^*(x)<+\infty$,
$\phi(x,\theta^*(x))\in V^c\cap{\bf X}$ and the infimum is attained. Moreover,
as mentioned above Definition \ref{JulyRem01}, $\phi(x,s)\in V^c$ for all $s\ge \theta^*(x)$. Therefore,
\begin{eqnarray*}
\int_{V}w(y)\eta(dy\times{\bf A})&=& \int_{V} w(x)\hat p(\bar\RR^0_+|x)\mu(dx\times\bar\RR^0_+\times{\bf A})\\
&&-\int_{V} w(x)\hat p([\theta^*(x),\infty]|x)\mu(dx\times\bar\RR^0_+\times{\bf A})\\
&&+\int_{V}\int_{\RR^0_+}\II\{\phi(x,\theta)\in V\}\int_{[0,\theta]}  \chi w(\phi(x,s))ds~\hat p(d\theta|x)\mu(dx\times\bar\RR^0_+\times{\bf A}).
\end{eqnarray*}
Recall, the measure $\eta(dy\times{\bf A})$ is finite and the function $\chi w(\phi(x,s))$ is integrable on $[0,\theta]$ with $\theta<\infty$.

After we apply the Tonelli Theorem \cite[Thm.11.28]{b4} to the last term, we obtain:
\begin{eqnarray*}
\int_{V}w(y)\eta(dy\times{\bf A})
&=& \int_{V} w(x)\mu(dx\times\bar\RR^0_+\times{\bf A})
-\int_{V} w(x)\hat p([\theta^*(x),\infty]|x)\mu(dx\times\bar\RR^0_+\times{\bf A})\\
&&+\int_{V}\int_{\RR^0_+}\int_{[s,\infty)} \II\{\phi(x,\theta)\in V\} \chi w(\phi(x,s))\hat p(d\theta|x)~ds~\mu(dx\times\bar\RR^0_+\times{\bf A})\\
&=& \int_{V} w(x)\mu(dx\times\bar\RR^0_+\times{\bf A})
-\int_{V} w(x)\hat p([\theta^*(x),\infty]|x)\mu(dx\times\bar\RR^0_+\times{\bf A})\\
&&+\int_{V}\int_{\RR^0_+}  \chi w(\phi(x,s))\II\{\phi(x,s)\in V\}\hat p([s,\theta^*(x))|x)~ds~\mu(dx\times\bar\RR^0_+\times{\bf A}).
\end{eqnarray*}
Note that
\begin{eqnarray*}&&\int_{V}\int_{\RR^0_+}  \chi w(\phi(x,s))\II\{\phi(x,s)\in V\}\hat p([s,\theta^*(x))|x)~ds~\mu(dx\times\bar\RR^0_+\times{\bf A})\\
&=&\int_{V}\int_{\RR^0_+}  \chi w(\phi(x,s))\hat p([s,\theta^*(x))|x)~ds~\mu(dx\times\bar\RR^0_+\times{\bf A})
\end{eqnarray*}
as $\chi w(\phi(x,s))=0$ for $\phi(x,s)\in V^c$. (See (\ref{ea1}), where, in our case, $V^c\subseteq D$ and $W(y)=0$ for all $y\in V^c$.)
Now
\begin{eqnarray*}
&&\int_{V}w(y)\eta(dy\times{\bf A})\\
&=& \int_{V} w(x)\mu(dx\times\bar\RR^0_+\times{\bf A})
-\int_{V} w(x)\hat p([\theta^*(x),\infty]|x)\mu(dx\times\bar\RR^0_+\times{\bf A})\\
&&+\int_{V} \chi w(y) \eta(dy\times\Box)-\int_{V}\int_{\RR^0_+} \chi w(\phi(x,s)) \hat p([\theta^*(x),\infty]|x) ~ds~ \mu(dx\times\bar\RR^0_+\times{\bf A}).
\end{eqnarray*}
All the integrals here are finite because, no matter whether $\theta^*(x)$ is finite or not,
\begin{eqnarray*}
\lim_{t\to\infty} w(\phi(x,t))=w(x)+\int_{\RR^0_+}\chi w(\phi(x,s))ds=0,~\forall~x\in V
\end{eqnarray*}
and thus
$$\int\limits_{V} w(x)\hat p([\theta^*(x),\infty]|x)\mu(dx\times\bar\RR^0_+\times{\bf A})+ \int\limits_{V}\int\limits_{\RR^0_+} \chi w(\phi(x,s)) \hat p([\theta^*(x),\infty]|x) ~ds~ \mu(dx\times\bar\RR^0_+\times{\bf A})=0.$$
This also leads to
\begin{eqnarray*}
\int_{V}w(y)\eta(dy\times{\bf A})
&=& \int_{V} w(x)\mu(dx\times\bar\RR^0_+\times{\bf A})+\int_{V} \chi w(y) \eta(dy\times\Box)\\
&=& w(x_0)+\int_{V\times\RR^0_+\times{\bf A}} w(l(\phi(y,\theta),a))\II\{\phi(y,\theta)\in V\}\mu(dy\times d\theta\times da)\\
&&+\int_{V} \chi w(y) \eta(dy\times\Box)
\end{eqnarray*}
by (\ref{e105}), and the required formula (\ref{e17}) follows from the definition (\ref{e11}). \hfill$\Box$

\section{Proof of Theorem \ref{t2} and Corollary \ref{corol1}}\label{sec6}

Below, we assume that Conditions \ref{ConstrainedPPZcondition05}, \ref{ConstrainedPPZcondition01}, \ref{co31} and \ref{con45} are satisfied.
The proofs will be based on a series of lemmas.

\begin{lemma}\label{l11}
Let $\pi=\{\pi_i\}_{i=1}^\infty$ be a reasonable Markov strategy as in Definition \ref{JulyRem01}, defined on $V$ by stochastic kernels $\pi_i(d\theta\times da|x)=p^i_T(d\theta|x) p^i_A(da|x,\theta)$. Suppose  $\eta$ is the corresponding aggregated occupation measure (\ref{e16p}) coming from the occupation measure $\mu^\pi$ as in (\ref{e101}).
Introduce the (partial) aggregated occupation measures
$$\eta^i(\Gamma_X\times\Gamma_A)\defi \eta^i(\Gamma_X\times(\Gamma_A\cap{\bf A}))+\eta^i(\Gamma_X\times\Box)\II\{\Box\in\Gamma_A\}$$
on $V\times {\bf A}_\Box$, defined recursively:
\begin{eqnarray*}
\eta^0(\Gamma_X\times\Gamma_A)&\equiv & 0;\\
\eta^{i+1}(\Gamma_X\times\Box)&=& \eta^i(\Gamma_X\times\Box)+\int_{V}\int_{\RR^0_+}\delta_{\phi(x,u)}(\Gamma_X)p^{i+1}_T([u,\infty]|x)du~\nu^i(dx),~\Gamma_X\in{\cal B}(V);\\
\eta^{i+1}(\Gamma_X\times\Gamma_A)&=& \eta^i(\Gamma_X\times\Gamma_A)+\int_{V}\int_{\RR^0_+}\delta_{\phi(x,\theta)}(\Gamma_X)p^{i+1}_A(\Gamma_A|x,\theta) p^{i+1}_T(d\theta|x)\nu^i(dx),\\
&&~~~~~\Gamma_X\in{\cal B}(V),~\Gamma_A\in{\cal B}({\bf A}),
\end{eqnarray*}
where $\nu^i(dx)=P^\pi_{x_0}(X_i\in dx)$ is the measure on $V$, $i\ge 0$.

Then $\eta^i\uparrow \eta$ on $V\times{\bf A}_\Box$ set-wise as $i\to\infty$. Every measure $\eta^i$ is normal.
\end{lemma}

\par\noindent\underline{Proof.} We will need the (partial) occupation measure on $V\times\bar\RR^0_+\times{\bf A}$
$$\mu^n(dx\times d\theta\times da)\defi E^\pi_{x_0}\left[\sum_{i=1}^n \II\{X_{i-1}\in dx,\Theta_i\in d\theta, A_i\in da\}\right],~~n=0,1,2,\ldots.$$
Clearly, $\mu^n\uparrow\mu^\pi$ on $V\times\bar\RR^0_+\times{\bf A}$ set-wise as $n\to\infty$. Therefore, according to the definition of the measure $\eta$, for each positive measurable function $C^g$ on $V$,
\begin{eqnarray*}
I^n&\defi& \int_{V\times\bar\RR^0_+\times{\bf A}}\left\{\int_{[0,\theta]} C^g(\phi(x,u))\II\{\phi(x,u)\in V\}du\right\} \mu^n(dx\times d\theta\times da)\\
&\uparrow & \int_{V\times\bar\RR^0_+\times{\bf A}}\left\{\int_{[0,\theta]} C^g(\phi(x,u))\II\{\phi(x,u)\in V\}du\right\} \mu^\pi(dx\times d\theta\times da)=\int_{V} C^g(y)\eta(dy\times\Box)
\end{eqnarray*}
and, for each positive measurable function $C^I$ on $V\times{\bf A}$,
\begin{eqnarray*}
J^n&\defi& \int_{V\times\bar\RR^0_+\times{\bf A}}\II\{\theta<+\infty\}\II\{\phi(x,\theta)\in V\} C^I(\phi(x,\theta),a) \mu^n(dx\times d\theta\times da)\\
&\uparrow & \int_{V\times\bar\RR^0_+\times{\bf A}}\II\{\theta<+\infty\}\II\{\phi(x,\theta)\in V\} C^I(\phi(x,\theta),a) \mu^\pi(dx\times d\theta\times da)\\
&=&\int_{V\times{\bf A}} C^I(y,a)\eta(dy\times da).
\end{eqnarray*}

We will prove by induction the following assertions:
$$I^n=\int_{V} C^g(y)\eta^n(dy\times\Box)~~\mbox{ and }~J^n=\int_{V\times{\bf A}} C^I(y,a)\eta^n(dy\times da).$$

If $n=0$, then $\mu^0=0$, $\eta^0=0$, $I^0=0$, and $J^0=0$.

Suppose the above assertions are valid for some $n\ge 0$. Then
\begin{eqnarray*}
I^{n+1}&=& I^n+\int_{V}\int_{\bar\RR^0_+}\int_{\bf A} \left\{\int_{[0,\theta]} C^g(\phi(x,u))\II\{\phi(x,u)\in V\} du\right\}p^{n+1}_A(da|x,\theta)p^{n+1}_T(d\theta|x)\nu^n(dx)\\
&&\mbox{and}\\
J^{n+1}&=& J^n+\int_{V}\int_{\bar\RR^0_+}\int_{\bf A} \II\{\theta<\infty\}\II\{\phi(x,\theta)\in V\} C^I(\phi(x,\theta),a)p^{n+1}_A(da|x,\theta)p^{n+1}_T(d\theta|x)\nu^n(dx)
\end{eqnarray*}
because on $V\times\bar\RR^0_+\times{\bf A}$ we have equality
$$\mu^{n+1}(dx\times d\theta\times da)=\mu^{n}(dx\times d\theta\times da)+p^{n+1}_A(da|x,\theta)p^{n+1}_T(d\theta|x)\nu^n(dx).$$
Recall that $Q(\{\Delta\}|x,f^*(x))=1$.
Using the Tonelli Theorem (see \cite[Thm.11.28]{b4}), we obtain:
\begin{eqnarray*}
&&\int_{V}\int_{\bar\RR^0_+}\left\{\int_{[0,\theta]} C^g(\phi(x,u))\II\{\phi(x,u)\in V\}du\right\} p^{n+1}_T(d\theta|x)\nu^n(dx)\\
&=&\int_{V}\int_{\RR^0_+}\int_{[u,\infty]}C^g(\phi(x,u))\II\{\phi(x,u)\in V\} p^{n+1}_T(d\theta|x)du~\nu^n(dx)\\
&=&\int_{V}\int_{\RR^0_+}C^g(\phi(x,u))\II\{\phi(x,u)\in V\} p^{n+1}_T([u,\infty]|x)du~\nu^n(dx)\\
&=&\int_{V} C^g(y)\left\{\int_{V}\int_{\RR^0_+} \delta_{\phi(x,u)}(dy)p^{n+1}_T([u,\infty]|x)du~\nu^n(dx)\right\},
\end{eqnarray*}
and, by induction and the definition of the measure $\eta^{n+1}(\Gamma_X\times\Box)$,
$$I^{n+1}=\int_{V} C^g(y)\eta^{n+1}(dy\times\Box).$$

Similarly,
\begin{eqnarray*}
&&\int_{V}\int_{\RR^0_+}\int_{\bf A} \II\{\phi(x,\theta)\in V\} C^I(\phi(x,\theta),a) p^{n+1}_A(da|x,\theta) p^{n+1}_T(d\theta|x)\nu^n(dx)\\
&=&\int_{V}\int_{\bf A} C^I(y,a)\left\{\int_{V}\int_{\RR^0_+} \delta_{\phi(x,\theta)}(dy)p^{n+1}_A(da|x,\theta) p^{n+1}_T(d\theta|x)\nu^n(dx)\right\},
\end{eqnarray*}
and, by induction and the definition of the measure $\eta^{n+1}(dx\times da)$ on $V\times \textbf{A}$,
\begin{eqnarray*}
J^{n+1}=\int_{V\times{\bf A}} C^I(y,a)\eta^{n+1}(dy\times da).
\end{eqnarray*}

Since, for all positive measurable functions $C^g$ on $V$ and $C^I$ on $V\times{\bf A}$,
\begin{eqnarray*}
&&\int_{V} C^g(y)\eta^n(dy\times\Box)\uparrow \int_{V} C^g(y)\eta(dy\times\Box)~\mbox{ and } \\
&&\int_{V\times{\bf A}} C^I(y,a)\eta^{n}(dy\times da)\uparrow \int_{V\times{\bf A}} C^I(y,a)\eta(dy\times da),
\end{eqnarray*}
we conclude that $\eta^n\uparrow\eta$ on $V\times{\bf A}_\Box$ set-wise as $n\to\infty$. The last assertion is obvious. \hfill $\Box$
\bigskip

\begin{lemma}\label{l9}
\begin{itemize}
\item
[(a)] Suppose  $\zeta$ is a finite measure on $V$, and the measure $\hat\zeta(d\tilde x^0)$ and the stochastic kernel $\check\zeta(dt|\tilde x^0)$ are as in Definition \ref{d8}. Then, for each bounded (or positive, or negative) measurable function $g$ on $V$,
\begin{eqnarray*}
\int_{V} g(y) \zeta(dy)=\int_{\tilde V}\int_{[0,\tilde\theta^*(\tilde x^0))} g(\phi((\tilde x^0,0),u))\check\zeta(du|\tilde x^0)\hat\zeta(d\tilde x^0)=\int_{\bf D} g(\phi((\tilde x^0,0),u))\check\zeta(d\tilde x^0\times du).
\end{eqnarray*}
\item[(b)] Suppose  $\zeta$ is a  normal measure on $V$, and the measure $\check\zeta(d\tilde x^0\times dt)$ is as in Definition \ref{d8}. Then, for each positive (or negative) measurable function $g$ on $V$,
\begin{eqnarray*}
\int_{V} g(y)\zeta(dy)=\int_{\bf D} g(\phi((\tilde x^0,0),u)\check\zeta(d\tilde x^0\times du).
\end{eqnarray*}
\item[(c)] Suppose $\zeta$ is a    normal (or finite) measure on the orbit
\begin{eqnarray*}
_{\tilde z}{\cal X}\cap V=\{\phi((\tilde z,0),t):~t\in[0,\tilde\theta^*(\tilde z))\}
\end{eqnarray*}
and
\begin{eqnarray*}
m(I)\defi\zeta(\{\phi((\tilde z,0),t):~t\in I\})
\end{eqnarray*}
is the $\sigma$-finite  (or finite)  measure on $[0,\tilde\theta^*(\tilde z))$. (The set $\{\phi((\tilde z,0),t):~t\in I\}$ is measurable because if $\tilde{z}\in\tilde{V}$, then $\phi((\tilde{z},0),\cdot)$ is a homeomorphism between $[0,\tilde{\theta}^\ast(\tilde{z}))$ and $_{\tilde z}{\cal X}\cap V$.) Then, for each positive  or negative  measurable function $g$ on $_{\tilde z}{\cal X}\cap V$,
\begin{eqnarray*}
\int_{_{\tilde z}{\cal X}\cap V} g(y)\zeta(dy)=\int_{[0,\tilde\theta^*(\tilde z))} g(\phi((\tilde z,0),t))m(dt).
\end{eqnarray*}
\end{itemize}
\end{lemma}

\par\noindent\underline{Proof.} (a) For the case of bounded functions $g$, it is sufficient to check the required formula for $g(y)=\II\{y\in Y\}$, where $Y\in{\cal B}(V)$ is an arbitrary set. According to the definition of the mappings $F$ and $F^{-1}$,
$$(\tilde x^0,u)\in F^{-1}(Y)\Longleftrightarrow F(\tilde x^0,u)\in Y\Longleftrightarrow (\tilde\phi(\tilde x^0,u),u)=\phi((\tilde x^0,0),u)\in Y.$$
Hence
\begin{eqnarray*}
\int_{V} g(y)\zeta(dy)&=& \zeta(Y)=\check{\zeta}(F^{-1}(Y))
= \int_{\tilde V\times\RR^0_+} \II\{(\tilde x^0,u)\in F^{-1}(Y)\} \check{\zeta}(d\tilde x^0\times du)\\
&=&\int_{\tilde V}\int_{\RR^0_+} \II\{\phi((\tilde x^0,0),u)\in Y\}\check\zeta(du|\tilde x^0)\hat\zeta(d\tilde x^0).
\end{eqnarray*}
Moreover, for $u\ge \tilde\theta^*(\tilde x^0)$, $\phi((\tilde x^0,0),u)\in V^c$ and thus $\phi((\tilde x^0,0),u)$ cannot belong to $Y$. The desired formula
$$\int_{V} g(y)\zeta(dy)=\int_{\tilde V}\int_{[0,\tilde\theta^*(\tilde x^0))} \II\{\phi((\tilde x^0,0),u)\in Y\}\check\zeta(du|\tilde x^0)\hat\zeta(d\tilde x^0)
=\int_{\bf D}\II\{\phi((\tilde x^0,0),u)\in Y\}\check\zeta(d\tilde x^0\times du)
$$
is proved.

For the case of positive functions $g$, one should apply the monotone convergence theorem to the sequence $g\wedge N\uparrow g$. Negative functions $g$ can be treated similarly.

(b) The required formula is justified after we represent the function $g$ as $g(x)=\sum_{t=1}^\infty g_t(x)$ with $g_t((\tilde x,u))=\II\{u\in[t-1,t)\} g((\tilde x,u))$ and use the statement (a) separately for  all $g_t$, where one can legitimately use the (finite) restriction of $\zeta$ to the set $\{x=(\tilde x,u)\in V:~t-1\le u<t\}$.

(c) Without loss of generality, we assume that $\zeta(_{\tilde z}{\cal X}\cap V)>0$. This implies $\tilde{z}\in\tilde{V}$ in particular.

If $\tilde\theta^*(\tilde z)<\infty$ then the measure $\zeta$ is finite and can be extended to $V$ by putting $\zeta(V\setminus ~_{\tilde z}{\cal X})\defi 0$. Now
\begin{eqnarray*}
\check\zeta(d\tilde x^0\times dt)&=&m(dt)\delta_{\tilde z}(d\tilde x^0);\\
\hat\zeta(d\tilde x^0)&=& \zeta(_{\tilde z}{\cal X}\cap V)\delta_{\tilde z}(d\tilde x^0);\\
\check\zeta(dt|\tilde x^0)&=&\left\{\begin{array}{ll}
m(dt)/\zeta(_{\tilde z}{\cal X}\cap V), & \mbox{ if } \tilde x^0=\tilde z;\\
\mbox{arbitrarily fixed probability measure}, & \mbox{ if } \tilde x^0\ne\tilde z,
\end{array}\right.
\end{eqnarray*}
and the required equality follows from Item (a). The same reasoning applies if $\tilde\theta^*(\tilde z)=\infty$ and the measure $\zeta$ is finite.

Suppose $\tilde\theta^*(\tilde z)=\infty$, so that $_{\tilde z}{\cal X}\cap V=~_{\tilde z}{\cal X}$, and the measure $\zeta$ is not finite, but normal.
It is sufficient to check the required formula for $g_t(y)=\II\{y\in Y_t\}$, where
\begin{eqnarray*}
Y_t=\{\phi((\tilde z,0),u):~u\in I_t\in{\cal B}([t-1,t))\},~~~t=1,2,\ldots.
\end{eqnarray*}
As mentioned in the statement of this lemma, the mapping $[G(u):=\phi((\tilde z,0),u)$ is a homeomorphism between $\RR^0_+$ and $_{\tilde z}{\cal X}$ (see Lemma \ref{la5}), and all different subsets $I_t\in{\cal B}([t-1,t))$ produce all possible subsets $Y_t\in{\cal B}(\{\phi((\tilde z,0),u):~u\in[t-1,t)\})$. Thus, for an arbitrary set $Y\in{\cal B}(_{\tilde z}{\cal X})$, we have $Y=\cup_{t=1}^\infty Y_t$ with
$$Y_t:=Y\cap\{\phi((\tilde z,0),u):~u\in[t-1,t)\}\in{\cal B}(\{\phi((\tilde z,0),u):~u\in[t-1,t)\}),$$
and the proof will be completed by applying the monotone convergence theorem.

Now
\begin{eqnarray*}
\int_{_{\tilde z}{\cal X}} g_t(y)\zeta(dy)&=&\zeta(Y_t);\\
\int_{\RR^0_+} g_t(\phi((\tilde z,0),u)) m(du) &=&\int_{\RR^0_+} \II\{\phi((\tilde z,0),u)\in Y_t\} m(du)=m(I_t),
\end{eqnarray*}
and $m(I_t)=\zeta(Y_t)$ by the definition of the measure $m$.
\hfill$\Box$
\bigskip

\begin{lemma}\label{l10}
Suppose an orbit
$$_{\tilde x}{\cal X}\cap V=\{\phi((\tilde x,0),t):~t\in[0,\tilde\theta^*(\tilde x))\}$$
is fixed and $p^*$ is a probability measure on $\bar\RR^0_+$ such that $p^*([\tilde\theta^*(\tilde x),\infty))=0$.

Then the measures  $\tilde\eta^*_\Box$ and $\tilde\eta^*_A$ on $~_{\tilde x}{\cal X}\cap V$, defined as
\begin{eqnarray*}
\eta^*_\Box(\Gamma) &\defi & \int_{\RR^0_+} \II\{\phi((\tilde x,0),u)\in\Gamma\} p^*([u,\infty]) du=\int_{\RR^0_+} \II\{\phi((\tilde x,0),u)\in\Gamma\} (1-p^*([0,u))du,\\
\eta^*_A(\Gamma) &\defi & \int_{\RR^0_+} \II\{\phi((\tilde x,0),u)\in\Gamma\}p^*(du),~~~~~~~~~~\Gamma\in{\cal B}(~_{\tilde x}{\cal X}\cap V),
\end{eqnarray*}
satisfy equation
\begin{equation}\label{e23prime}
0=w((\tilde x,0))+\int_{~_{\tilde x}{\cal X}\cap V}\chi w(x)\eta^*_\Box(dx)-\int_{~_{\tilde x}{\cal X}\cap V} w(x)\eta^*_A(dx)
\end{equation}
for all functions $w\in{\bf W}$. The measure $\eta^*_A$ is finite, and the measure $\eta^*_\Box$ is normal on that orbit.
\end{lemma}

\par\noindent\underline{Proof.}
The properties of the measures $\eta^*_A$ and $\eta^*_\Box$  formulated  in the last sentence of this lemma are obvious, c.f. the reasoning in the proof of Lemma \ref{l101}(a).

Now let $w\in\textbf{W}$ be fixed. We verify the rest of the statement of this lemma by distinguishing the following two cases.

(i) Suppose that $u^*\defi\inf\{u\in\bar\RR^0_+:~p^*([0,u])=1\}\ge \tilde\theta^*(\tilde x)$.
The expression
\begin{eqnarray*}
I\defi w((\tilde x,0))+\int_{~_{\tilde x}{\cal X}\cap V}\chi w(x)\eta^*_\Box(dx)-\int_{~_{\tilde x}{\cal X}\cap V} w(x)\eta^*_A(dx)
\end{eqnarray*}
is well defined because the measure $\eta^*_\Box$ is normal, the integral $\int_{~_{\tilde x}{\cal X}\cap V}\chi w(x)\eta^*_\Box(dx)$ is positive or negative,
the function $w$ is bounded and the measure $\eta^*_A$ is finite. According to Lemma \ref{l9}(c),
\begin{eqnarray*}
I &= &  w((\tilde x,0))+\int_{[0,\tilde\theta^*(\tilde x))} \chi w(\phi((\tilde x,0),t)) [1- p^*([0,t))]~ dt-\int_{[0,\tilde\theta^*(\tilde x))} w(\phi((\tilde x,0),t)) p^*(dt)\\
&=& -\left[\int_{[0,\tilde\theta^*(\tilde x))} \chi w(\phi((\tilde x,0),t)) p^*([0,t))~ dt+\int_{[0,\tilde\theta^*(\tilde x))} w(\phi((\tilde x,0),t)) p^*(dt) \right].
\end{eqnarray*}
The last equality is by Lemma \ref{l1} and Definition \ref{d3} of the space $\bf W$:
$$w((\tilde x,0))+\lim_{T\to\tilde\theta^*(\tilde x)} \int_{[0,T]} \chi w(\phi((\tilde x,0),t))dt=\lim_{T\to\tilde\theta^*(\tilde x)} w(\phi((\tilde x,0),T))=0.$$
We apply the Tonelli Theorem \cite[Thm.11.28]{b4} to the first integral in the square brackets  and again use Lemma \ref{l1}:
\begin{eqnarray*}
\int_{[0,\tilde\theta^*(\tilde x))} \chi w(\phi((\tilde x,0),t))\int_{[0,t)} p^*(du)~dt&=&\int_{[0,\tilde\theta^*(\tilde x))} \int_{(u,\tilde\theta^*(\tilde x))} \chi w(\phi((\tilde x,0),t))dt~ p^*(du)\\
&=& \int_{[0,\tilde\theta^*(\tilde x))} [-w(\phi((\tilde x,0),u))] p^*(du).
\end{eqnarray*}
Thus $I=0$.

(ii) Suppose that $u^*\defi\inf\{u\in\bar\RR^0_+:~p^*([0,u])=1\}< \tilde\theta^*(\tilde x)$. Since measures $\tilde\eta^\ast_A$ and $\tilde\eta^\ast_\Box$ both equal zero on the set $\{\phi((\tilde x,0),t):~t>u^*\}$, it is sufficient to show that
$$I\defi w((\tilde x,0))+\int_{{\cal X}_0^{u^*}}\chi w(x)\eta^*_\Box(dx)-\int_{{\cal X}_0^{u^*}} w(x)\eta^*_A(dx)=0,$$
where
\begin{eqnarray*}
{\cal X}_0^{u^*}\defi\{\phi((\tilde x,0),t):~0\le t\le u^*\}.
\end{eqnarray*}
This expression is well defined because the measure $\eta^*_\Box$ is normal, the integral $\int_{{\cal X}_0^{u^*}}\chi w(x)\eta^*_\Box(dx)$ is positive or negative, the function $w$ is bounded and the measure $\eta^*_A$ is finite. The measure $\eta^*_\Box$ is non-atomic, and the first integral can be calculated over
\begin{eqnarray*}
{\cal X}_{0}^{u^*-}\defi\{\phi((\tilde x,0),t):~0\le t< u^*\},
\end{eqnarray*}
so that, by Lemma \ref{l9}(c),
\begin{eqnarray*}
I&=&w((\tilde x,0))+\int_{[0,u^*)} \chi w(\phi((\tilde x,0),t)) [1-p^*([0,t))]dt\\
&&-\int_{[0,u^*)} w(\phi((\tilde x,0),t)) p^*(dt)-w(\phi((\tilde x,0),u^*))[1-p^*([0,u^*))].
\end{eqnarray*}
In the last term, $[1-p^*([0,u^*))]=p^*(\{u^*\})$. Since
$$w((\tilde x,0))+\int_{[0,u^*)} \chi w(\phi((\tilde x,0),t)) dt-w(\phi((\tilde x,0),u^*))=0$$
(see Lemma \ref{l1}), after we subtract  this equality from $I$, we obtain
$$I=-\int_{[0,u^*)} \chi w(\phi((\tilde x,0),t)) p^*([0,t))dt-\int_{[0,u^*)} w(\phi((\tilde x,0),t)) p^*(dt)+w(\phi((\tilde x,0),u^*))p^*([0,u^*)).$$
Finally, apply the Tonelli Theorem (see \cite[Thm.11.28]{b4}) to the first term and again use Lemma \ref{l1}:
\begin{eqnarray*}
\int_{[0,u^*)}\int_{[0,t)} \chi w(\phi((\tilde x,0),t))~p^*(du)~dt &=&  \int_{[0,u^*)}\int_{(u,u^*)} \chi w(\phi((\tilde x,0),t))~dt~p^*(du)\\
&=& \int_{[0,u^*)} [w(\phi((\tilde x,0),u^*))-w(\phi((\tilde x,0),u))] p^*(du)\\
&=& w(\phi((\tilde x,0),u^*))p^*([0,u^*))-\int_{[0,u^*)} w(\phi((\tilde x,0),u))~ p^*(du).
\end{eqnarray*}
Therefore, $I=0$.

The proof is completed.
\hfill$\Box$

\begin{lemma}\label{l8}
Suppose $\nu$ is a finite measure on $V$ such that $\nu(V\cap (\tilde V\times\{t:~t>0\}))=0$,
$\tilde\eta$ is a finite measure on $V\times{\bf A}$,  $\tilde\eta_\Box$  is a normal  measure on $V$ and $\tilde\eta_A$ is a finite measure on $V$ which satisfy equation
\begin{equation}\label{e20}
0=\int_{V} w(x)\nu(dx)+\int_{V}\chi w(x)\tilde\eta_\Box(dx)-\int_{V} w(x)\tilde\eta_A(dx)+\int_{V\times{\bf A}} w(l(x,a)) \tilde\eta(dx\times da)
\end{equation}
for all functions $w\in{\bf W}$.
Then there is a stochastic kernel $\tilde p(dt|x)$ on $\bar\RR^0_+$ given $V$ such that, for $\theta^*$ given by
\begin{equation}\label{e45}
\theta^*(x)\defi \inf\{\theta\in\RR^0_+:~\phi(x,\theta)\in V^c\},
\end{equation}
$\tilde p([\theta^*(x),\infty)|x)=0$ for all $x\in V$
and the measures
\begin{eqnarray*}
\tilde\eta'_A(\Gamma)&\defi& \int_{V} \int_{\RR^0_+} \II\{\phi(x,u)\in\Gamma\} \tilde p(du|x) \nu(dx)\\
\mbox{and }~ \tilde\eta'_\Box(\Gamma)&\defi & \int_{V} \int_{\RR^0_+} \II\{\phi(x,u)\in\Gamma\} \tilde p([u,\infty]|x)du~ \nu(dx),~~\Gamma\in{\cal B}(V)
\end{eqnarray*}
satisfy equation
\begin{equation}\label{e13}
0=\int_{V} w(x)\nu(dx)+\int_{V}\chi w(x)\tilde\eta'_\Box(dx)-\int_{V} w(x)\tilde\eta'_A(dx)
\end{equation}
for all functions $w\in{\bf W}$.
Moreover, the set functions $\tilde\eta_\Box(\Gamma)-\tilde\eta'_\Box(\Gamma)$ and $\tilde\eta_A(\Gamma)-\tilde\eta'_A(\Gamma)$ on ${\cal B}(V)$ are again  normal  and finite measures, correspondingly.
\end{lemma}

\par\noindent\underline{Proof.} (i) Firstly, we introduce several functions, measures and sets , describe their properties and define explicitly the stochastic kernel $\tilde p$.

The necessary properties of the function $\theta^*$ were established during the proof of Theorem \ref{t1}. Note that, for each $\Gamma\in{\cal B}(V)$, the function $\II\{\phi(x,u)\in\Gamma\}$ is measurable since the flow $\phi$ is continuous. Below, $\hat\nu(\tilde\Gamma)\defi \nu(\tilde\Gamma\times\{0\})$ for $\tilde\Gamma\in{\cal B}(\tilde V)$.

In accordance with Definition \ref{d8}, we introduce the finite measure $\hat\eta_A(d\tilde x^0)$ and stochastic kernel
$\check\eta_A(dt|\tilde x^0)$ coming from $\tilde\eta_A(dx)$.
Next, introduce the finite measure
$$K\defi \hat\nu+\hat\eta_A$$
on $\tilde V$ and the Radon-Nikodym derivatives
$$n(\tilde x^0)\defi\frac{d\hat\nu}{dK}(\tilde x^0),~\mbox{ and }~a(\tilde x^0)\defi \frac{d\hat\eta_A}{dK}(\tilde x^0).$$
Below, we fix one specific version of the derivative $n$ and of the derivative $a$. On the set
$$\tilde{\bf V}_\nu\defi \{\tilde x^0\in \tilde V:~n(\tilde x^0)>0\},$$
we have
$$\hat\eta_A(\tilde\Gamma)= \int_{\tilde\Gamma} a(\tilde x^0)K(d\tilde x^0)=\int_{\tilde\Gamma} \frac{a(\tilde x^0)}{n(\tilde x^0)} \hat\nu(d\tilde x^0)$$
for all $\tilde\Gamma\in{\cal B}(\tilde{\bf V}_\nu)$. 
See Figure \ref{fig3}. Note that $\hat{\nu}(\tilde V\setminus\tilde{\bf V}_\nu)=0$.

Since the function $\II\{u\le t\}$ of $(u,t)$ is measurable, the integral $\int_{\RR^0_+} \II\{u\le t\}\check\eta_A(du|\tilde x^0)$ is a measurable function of $(\tilde x^0,t)$ (see \cite[Prop.7.29]{Bertsekas:1978}), and hence the function
\begin{eqnarray*}
G(\tilde x^0,t)\defi\check\eta_A([0,t]|\tilde x^0)\frac{a(\tilde x^0)}{n(\tilde x^0)}=\int_{\RR^0_+} \II\{u\le t\} \check\eta_A(du|\tilde x^0) \frac{a(\tilde x^0)}{n(\tilde x^0)},~~~\tilde x^0\in\tilde{\bf V}_\nu,~t\in\RR^0_+
\end{eqnarray*}
is measurable. For all $\tilde x^0\in\tilde{\bf V}_\nu$, the function $G(\tilde x^0,\cdot)$ clearly increases and is right-continuous: it is constant for $t\ge\tilde\theta^*(\tilde x^0)$ and, if $t_i\downarrow t\in[0,\tilde\theta^*(\tilde x^0))$ then $\check\eta_A([0,t_i]|\tilde x^0)\downarrow \check\eta_A([0,t]|\tilde x^0)$.

\begin{figure}[!htb]
	\centering
	\includegraphics[scale=0.3]{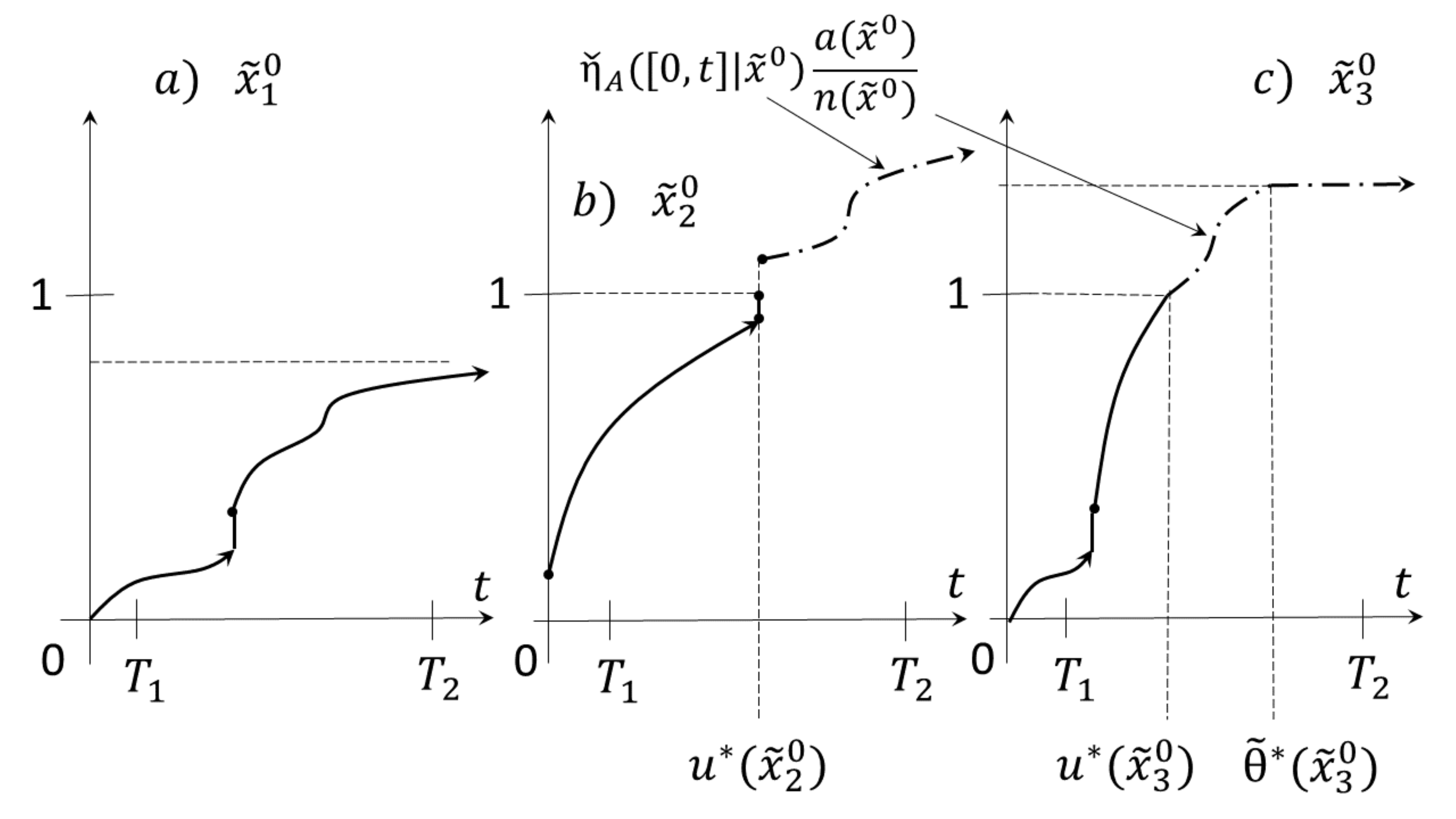}
	\caption{Graphs of the function $\tilde p([0,t]|(\tilde x^0,0))$, see also Figure \ref{fig5}.
In case a), $\tilde p([0,t]|(\tilde x^0_1,0))=\check\eta_A([0,t]|\tilde x^0_1)\frac{a(\tilde x^0_1)}{n(\tilde x^0_1)}$ for all $t\in\RR^0_+$, $u^*(\tilde x^0_1)=\tilde\theta^*(\tilde x^0_1)=\infty$ and $\tilde p(\RR^0_+|(\tilde x^0_1,0))<1$. In case b), $0<u^*(\tilde x^0_2)<\tilde\theta^*(\tilde x^0_2)$, $\check\eta_A(\{u^*(\tilde x^0_2)\}|\tilde x^0_2)>0$. In case c), $0<u^*(\tilde x^0_3)<\tilde\theta^*(\tilde x^0_3)<\infty$, $\check\eta_A(\{u^*(\tilde x^0_3)\}|\tilde x^0_3)=0$.  }
	\label{fig7}
\end{figure}

Let us introduce the function
\begin{eqnarray*}
u^*(\tilde x^0)\defi\inf\{t\in\mathbb{R}_+^0:~G(\tilde x^0,t)\ge 1\}\in\bar\RR^0_+,~~~\tilde x^0\in\tilde{\bf V}_\nu.
\end{eqnarray*}
When $u^*(\tilde x^0)<\infty$, this infimum is  attained because  the function $G(\tilde x^0,\cdot)$ is right-continuous; and $G(\tilde x^0,u^*(\tilde x^0)-)\le 1$. To show that the function $u^*(\cdot)$ is measurable, note that the function
$$f(\tilde x^0,t):=\infty\times\II\{G(\tilde x^0,t)<1\}+t\times\II\{G(\tilde x^0,t)\ge 1\}$$
is measurable and the function $t\to f(\tilde x^0,t)$ is lower semicontinuous for each $\tilde x^0\in\tilde{\bf V}_\nu$. Now, the function $u^*(\tilde x^0)=\inf_{t\in\bar\RR^0_+} f(\tilde x^0,t)$ is measurable by \cite[Thm.2]{Himmelberg:1976}; see also Corollary 1 and Remark 1 of \cite{Brown:1973}.
Note also that if $u^*(\tilde x^0)>\tilde\theta^*(\tilde x^0)$, then $u^*(\tilde x^0)=\infty$. Figure \ref{fig7} can serve as an illustration.

For $\tilde x^0\in\tilde{\bf V}_\nu$, we put
\begin{eqnarray*}
\tilde p(I|(\tilde x^0,0)) &\defi & \check\eta_A(I\cap[0,u^*(\tilde x^0)\wedge\tilde\theta^*(\tilde x^0))|\tilde x^0)\frac{a(\tilde x^0)}{n(\tilde x^0)}\\
&&+\II\{u^*(\tilde x^0)<\tilde\theta^*(\tilde x^0)\}\II\{u^*(\tilde x^0)\in I\}\left[1-\check\eta_A([0,u^*(\tilde x^0))|\tilde x^0)\frac{a(\tilde x^0)}{n(\tilde x^0)}\right]\\
&&\mbox{for all }~I\in{\cal B}(\RR^0_+),\\
\mbox{and}~\tilde p(\{\infty\}|(\tilde x^0,0))&\defi& 1-\tilde p(\RR^0_+|(\tilde x^0,0)).
\end{eqnarray*}
For all other points $x\in V$, we put $\tilde p(\{\infty\}|x)=1$ and $\tilde p(I|x)\equiv 0$ for $I\in{\cal B}(\RR^0_+)$. Clearly, $\tilde p([\theta^*(x),\infty)|x)=0$ for  all $x\in V$. The possible shapes of the distribution function $\tilde p([0,t]|(\tilde x^0,0))$ are shown on Figure \ref{fig7}.

(ii) Let us prove that equation (\ref{e13}) holds. Since $\nu(\tilde V\times\{t:~t>0\})=0$,
\begin{eqnarray*}
\tilde\eta'_A(\Gamma)=\int_{\tilde V}\eta^*_A(\Gamma|(\tilde x^0,0))\hat\nu(d\tilde x^0);~~~
\mbox{and }~\tilde\eta'_\Box(\Gamma) = \int_{\tilde V}\eta^*_\Box(\Gamma|(\tilde x^0,0))\hat\nu(d\tilde x^0)
\end{eqnarray*}
for all $\Gamma\in{\cal B}(V)$,
where
\begin{eqnarray*}
\eta^*_A(\Gamma|(\tilde x^0,0))&\defi&  \int_{\RR^0_+} \II\{\phi((\tilde x^0,0),u)\in\Gamma\} \tilde p(du|(\tilde x^0,0));\\
\eta^*_\Box(\Gamma|(\tilde x^0,0))&\defi& \int_{\RR^0_+} \II\{\phi((\tilde x^0,0),u)\in\Gamma\} \tilde p([u,\infty]|(\tilde x^0,0))du.
\end{eqnarray*}
The introduced measures $\eta^*_A$ and $\eta^*_\Box$ are concentrated on $~_{\tilde x^0}{\cal X}\cap V$ for each $\tilde x^0\in\tilde V$. By the way, $\eta^*_A$ and $\eta^*_\Box$ are measurable kernels because the flow $\phi$ is continuous and $\tilde p$ is a (measurable) stochastic kernel. Now
\begin{eqnarray*}
&&\int_{V} w(x)\nu(dx)+\int_{V}\chi w(x)\tilde\eta'_\Box(dx)-\int_{V} w(x)\tilde\eta'_A(dx)\\
&=& \int_{\tilde V} w((\tilde x^0,0))\hat\nu(d\tilde x^0)+\int_{\tilde V}\int_{_{\tilde x^0}{\cal X}\cap V} \chi w(x)\eta^*_\Box(dx|(\tilde x^0,0))\hat\nu(d\tilde x^0)\\
&&-\int_{\tilde V}\int_{_{\tilde x^0}{\cal X}\cap V} w(x)\eta^*_A(dx|(\tilde x^0,0))\hat\nu(d\tilde x^0)\\
&=& \int_{\tilde V}\left[w((\tilde x^0,0))+\int_{_{\tilde x^0}{\cal X}\cap V} \chi w(x) \eta^*_\Box(dx|(\tilde x^0,0))-\int_{_{\tilde x^0}{\cal X}\cap V} w(x)\eta^*_A(dz|(\tilde x^0,0))\right] \hat\nu(d\tilde x^0).
\end{eqnarray*}
The re-arrangement is legal because the function $w$ is bounded,  the function $\chi w$ is positive (or negative), the measure $\tilde\eta'_\Box$ is normal, and the measures $\hat\nu$ and $\eta^*_A(dx|(\tilde x^0,0))$ are finite (for all $\tilde x^0\in \tilde{V}$). Equation (\ref{e13}) follows from Lemma \ref{l10}.

(iii) Let us show that $\tilde\eta_A -\tilde\eta'_A $ is a finite measure. In case $\tilde x^0\in\tilde{\bf V}_\nu$ and $u^*(\tilde x^0)<\tilde\theta^*(\tilde x^0)$,
$$\tilde p(\{u^*(\tilde x^0)\}|(\tilde x^0,0))=1-\check\eta_A([0,u^*(\tilde x^0))|\tilde x^0)\frac{a(\tilde x^0)}{n(\tilde x^0)}\le \check\eta_A(\{u^*(\tilde x^0)\}|\tilde x^0)\frac{a(\tilde x^0)}{n(\tilde x^0)}$$
because
$$\check\eta_A([0,u^*(\tilde x^0)]|\tilde x^0)\frac{a(\tilde x^0)}{n(\tilde x^0)}=G(\tilde x^0,u^*(\tilde x^0))\ge 1.$$
Therefore, whether $u^*(\tilde x^0)<\tilde\theta^*(\tilde x^0)$ or $u^*(\tilde x^0)\ge\tilde\theta^*(\tilde x^0)$,
\begin{eqnarray*}
\tilde p(I|(\tilde x^0,0))\le\check\eta_A(I|\tilde x^0)\frac{a(\tilde x^0)}{n(\tilde x^0)}
\end{eqnarray*}
for all $I\in{\cal B}(\RR^0_+)$ and for  all $\tilde x^0\in\tilde{\bf V}_\nu$. Now, for each measurable subset $\Gamma\subset V$,
\begin{eqnarray*}
\tilde\eta'_A(\Gamma)&=& \int_{\tilde{\bf V}_\nu}\int_{\RR^0_+} \II\{\phi((\tilde x^0,0),u)\in\Gamma\}\tilde p(du|(\tilde x^0,0))\hat\nu(d\tilde x^0)\\
&=&\int_{\tilde{\bf V}_\nu}\int_{[0,\tilde\theta^*(\tilde x^0))} \II\{\phi((\tilde x^0,0),u)\in\Gamma\}\tilde p(du|(\tilde x^0,0))\hat\nu(d\tilde x^0)\\
&\le & \int_{\tilde{\bf V}_\nu}\int_{[0,\tilde\theta^*(\tilde x^0))}\II\{\phi((\tilde x^0,0),u)\in\Gamma\} \check\eta_A(du|\tilde x^0)\frac{a(\tilde x^0)}{n(\tilde x^0)} \hat\nu(d\tilde x^0)\\
&=&\int_{\tilde{\bf V}_\nu}\int_{[0,\tilde\theta^*(\tilde x^0))}\II\{\phi((\tilde x^0,0),u)\in\Gamma\} \check\eta_A(du|\tilde x^0) \hat\eta_A(d\tilde x^0)\\
&\le &\int_{\tilde V}\int_{[0,\tilde\theta^*(\tilde x^0))}\II\{\phi((\tilde x^0,0),u)\in\Gamma\} \check\eta_A(du|\tilde x^0) \hat\eta_A(d\tilde x^0)\\
&=& \int_{V} \II\{y\in\Gamma\}\tilde\eta_A(dy)=\tilde\eta_A(\Gamma).
\end{eqnarray*}
The last but one equality is by Lemma \ref{l9}(a). Hence, $\tilde\eta_A-\tilde\eta'_A$ is a finite measure.

(iv) Let us show that $\tilde\eta_\Box\ge \tilde\eta'_\Box$ set-wise. Recall that the measure $\tilde\eta'_\Box$ is normal. It is convenient to consider, with some abuse of notations, the images $\check\eta_\Box$, $\check\eta'_\Box$ and $\check\eta$ of the measures $\tilde\eta_\Box$, $\tilde\eta'_\Box$ and $\tilde\eta(\cdot\times{\bf A})$ as in Definition \ref{d8}. Recall that $\tilde\eta_\Box\ge \tilde\eta'_\Box\Leftrightarrow \check\eta_\Box\ge \check\eta'_\Box$.
Now, according to Lemma \ref{l9}(a,b), equation (\ref{e20}) takes the form:
\begin{eqnarray}
0&=& \int_{\tilde V} w((\tilde x^0,0))\hat\nu(d\tilde x^0)+\int_{\bf D} \chi w(\phi((\tilde x^0,0),u))\check\eta_\Box(d\tilde x^0\times du)\label{e21}\\
&&-\int_{\tilde V}\int_{[0,\tilde\theta^*(\tilde x^0))} w(\phi((\tilde x^0,0),u)) \check\eta_A(du|\tilde x^0)\hat\eta_A(d\tilde x^0)
 +\int_{\bf D} w^A(\phi((\tilde x^0,0),u)) \check\eta(d\tilde x^0\times du),\nonumber
\end{eqnarray}
where
\begin{equation}\label{e22}
w^A(y)\defi \int_{\bf A} w(l(y,a))\tilde\eta^{A}(da|y)
\end{equation}
and the stochastic kernel $\tilde\eta^{A}(da|y)$ comes from the decomposition
\begin{eqnarray*}
\tilde\eta(dy\times da)=\tilde\eta^{A}(da|y)\tilde\eta(dy\times{\bf A}).
\end{eqnarray*}

\begin{figure}[!htb]
	\centering
	\includegraphics[scale=0.3]{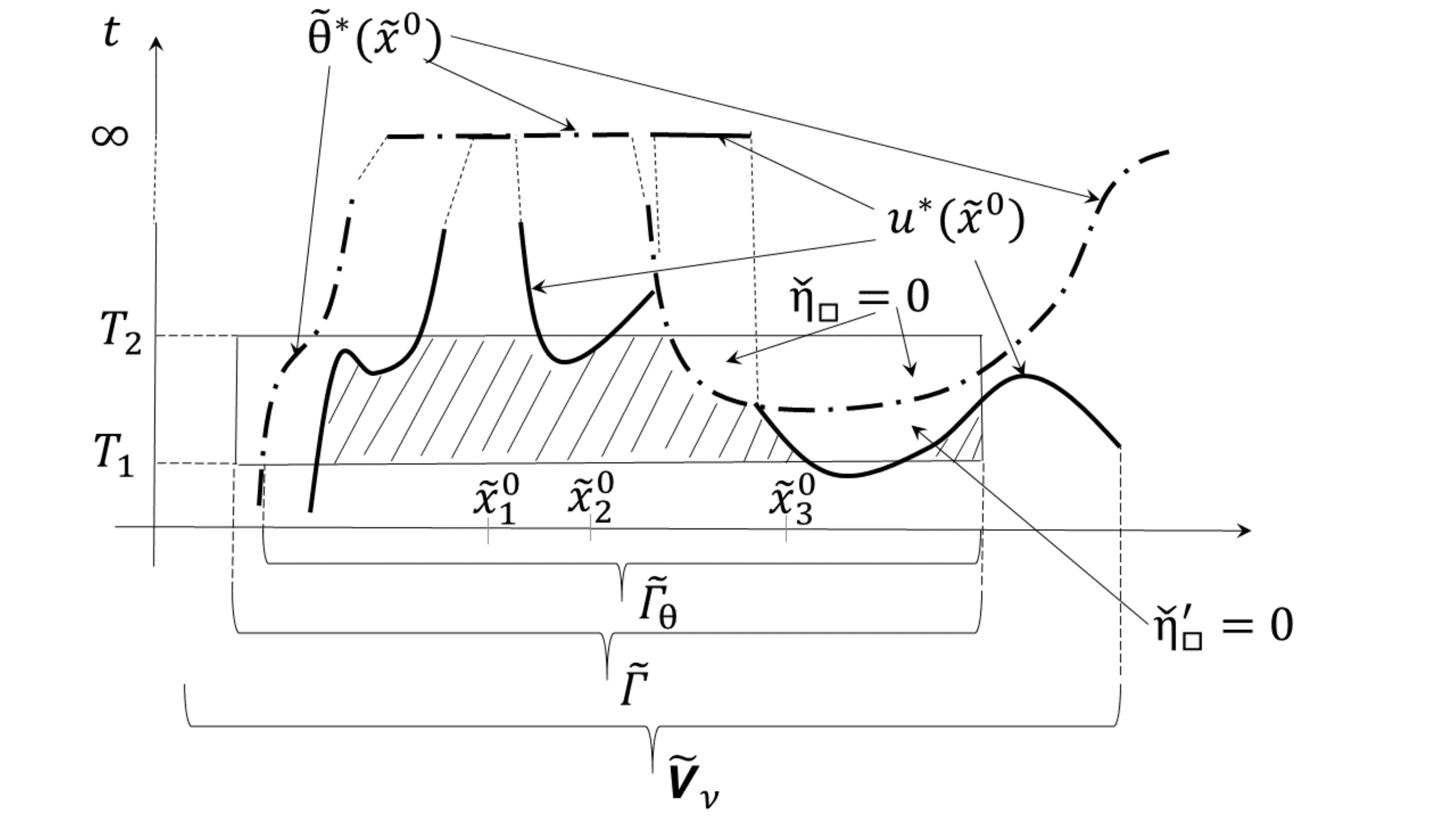}
	\caption{Space ${\bf D}=\{(\tilde x^0,t):~\tilde\phi(\tilde x^0,t)\in\tilde V\}$ and ``rectangle'' $Y_{T_1,T_2,\tilde\Gamma}$. The points $\tilde{x}_1^0,\tilde{x}_2^0,\tilde{x}_3^0$ belong to $\tilde{\bf V}_\nu\subset \tilde{V}.$ The dashed area is the part of $Y_{T_1,T_2,\tilde\Gamma}$ where $\check\eta'_\Box$ might be positive.}
	\label{fig5}
\end{figure}

According to \cite[V.1;Thm.1.5.6]{b5}, it suffices to show that the value of the measure $\check\eta_\Box$ is greater or equal to the value of $\check\eta'_\Box$  on each set of the form
\begin{eqnarray*}
Y_{T_1,T_2,\tilde\Gamma}\defi\{(\tilde x^0,u):~ \tilde x^0\in\tilde\Gamma,~T_1\le u<T_2\wedge\tilde\theta^*(\tilde x^0)\},~~~\tilde\Gamma\in{\cal B}(\tilde V),~0\le T_1<T_2<\infty.
\end{eqnarray*}
See Figure \ref{fig5} and also Figure \ref{fig7} for illustration.

Note that, in case $\tilde\Gamma\subset \tilde V\setminus\tilde{\bf V}_\nu$, since $\nu((\tilde V\setminus\tilde{\bf V}_\nu)\times\{0\})=0$,
$\check\eta'_\Box(Y_{T_1,T_2,\tilde\Gamma})=0$ and hence $\check\eta_\Box(Y_{T_1,T_2,\tilde\Gamma})-\check\eta'_\Box(Y_{T_1,T_2,\tilde\Gamma})\ge 0$ for all $T_1,T_2$. Therefore, below in this proof, we assume that
$\tilde\Gamma\subset\tilde{\bf V}_\nu$.

\begin{figure}[!htb]
	\centering
	\includegraphics[scale=0.2]{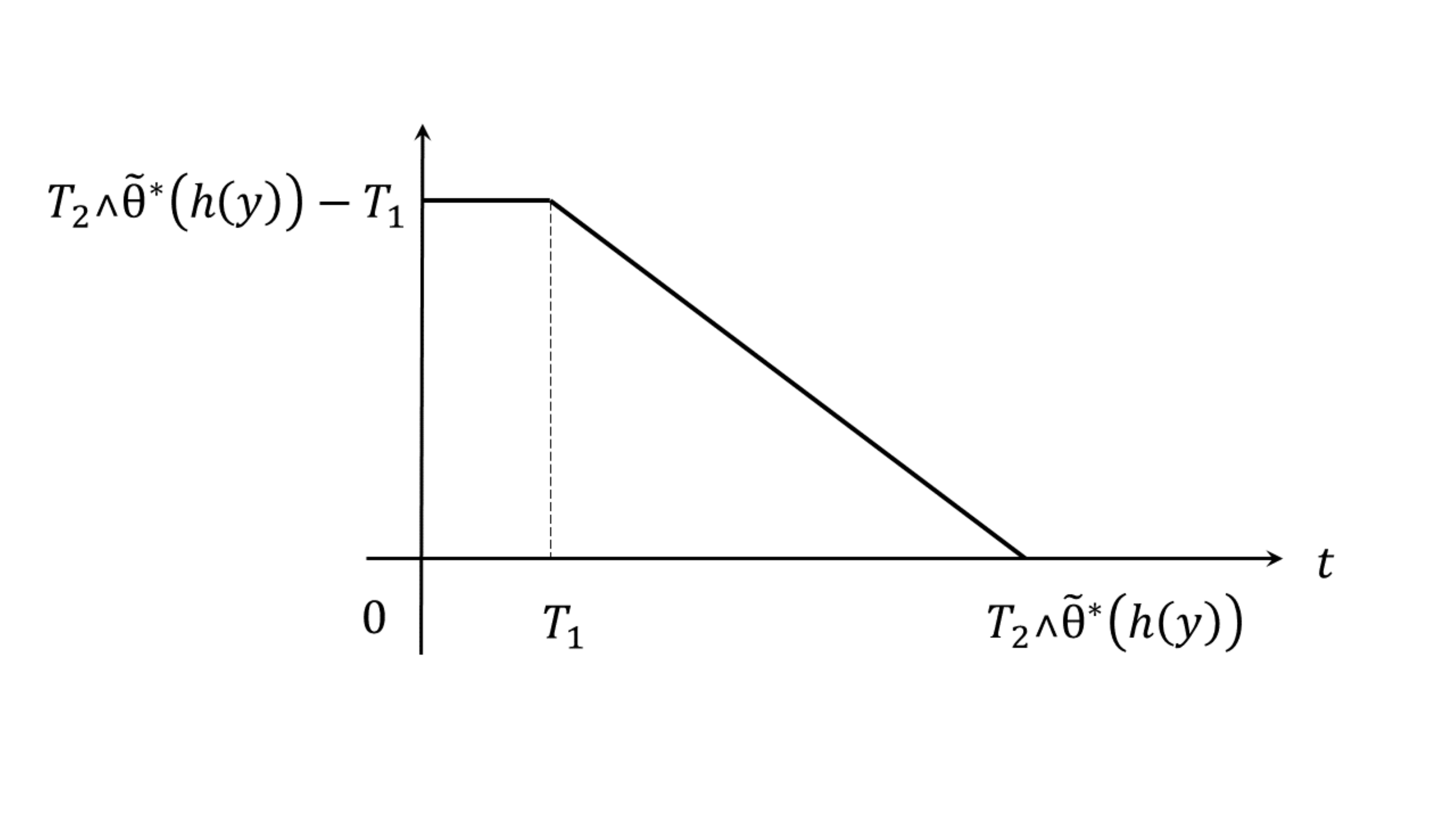}
	\caption{Graph of the function $w_{T_1,T_2,\tilde\Gamma}((\tilde y,t))=w_{T_1,T_2,\tilde\Gamma}((\tilde\phi(\tilde x^0,t),t))$
for a fixed value of $h(y)=\tilde x^0\in\tilde\Gamma$ and $\tilde\theta^*(h(y))>T_1$.}
	\label{fig6}
\end{figure}

To use equality (\ref{e21}) for calculating $\check\eta_\Box(Y_{T_1,T_2,\tilde\Gamma})$, we put
\begin{eqnarray*}
\chi w_{T_1,T_2,\tilde\Gamma}(\phi((\tilde x^0,0),u))\defi -\II\{(\tilde x^0,u)\in Y_{T_1,T_2,\tilde\Gamma}\},
\end{eqnarray*}
and  consider the following positive function decreasing along the flow:
\begin{eqnarray}\label{e27}
&&w_{T_1,T_2,\tilde\Gamma}(y)=w_{T_1,T_2,\tilde\Gamma}((\tilde y,t))\\
&\defi& \II\{h(y)\in\tilde\Gamma\}\times \left\{\begin{array}{ll}
T_2\wedge\tilde\theta^*(h(y))-T_1\wedge\tilde\theta^*(h(y)), & \mbox{ if } 0\le t\le T_1;\\
T_2\wedge\tilde\theta^*(h(y))-t, & \mbox{ if } T_1<t\le T_2\wedge\tilde\theta^*(h(y));\\
0,& \mbox { if } t>T_2\wedge\tilde\theta^*(h(y)).
\end{array}\right. \nonumber
\end{eqnarray}
See Figure \ref{fig6}. The function $h$ was introduced in Definition \ref{d7}.
Clearly,  $w_{T_1,T_2,\tilde\Gamma}\in{\bf W}$ for all $0\le T_1<T_2<\infty$, $\tilde\Gamma\in{\cal B}(\tilde V)$, and
\begin{eqnarray*}
w_{T_1,T_2,\tilde\Gamma}((\tilde x^0,0))=\II\{\tilde x^0\in\tilde\Gamma\}\left(T_2\wedge \tilde\theta^*(\tilde x^0)-T_1\wedge \tilde\theta^*(\tilde x^0)\right);
\end{eqnarray*}
\begin{eqnarray*}w_{T_1,T_2,\tilde\Gamma}(\phi((\tilde x^0,0),u))=\II\{\tilde x^0\in\tilde\Gamma\}\times\left\{\begin{array}{ll}
T_2\wedge\tilde\theta^*(\tilde x^0)-T_1\wedge\tilde\theta^*(\tilde x^0), & \mbox{ if } u\le T_1;\\
T_2\wedge\tilde\theta^*(\tilde x^0)-u, & \mbox{ if } T_1<u\le T_2\wedge\tilde\theta^*(\tilde x^0);\\
0,& \mbox { if } u>T_2\wedge\tilde\theta^*(\tilde x^0).
\end{array}\right.
\end{eqnarray*}
The expression (\ref{e22}) takes the form
\begin{eqnarray*}
w^A_{T_1,T_2,\tilde\Gamma}(y)&\defi&  \int_{\bf A} w_{T_1,T_2,\tilde\Gamma}( l(y,a))\tilde\eta^{A}(da|y).
\end{eqnarray*}
From equality (\ref{e21}), using the expression $\hat\eta_A(d\tilde x^0)=\frac{a(\tilde x^0)}{n(\tilde x^0)}\hat\nu(d\tilde x^0)$, we have for
\begin{eqnarray*}
\tilde\Gamma_\theta\defi \tilde\Gamma\cap\{\tilde x^0:~\tilde\theta^*(\tilde x^0)\ge T_1\}:
\end{eqnarray*}
\begin{eqnarray*}
\check\eta_\Box(Y_{T_1,T_2,\tilde\Gamma})&=& \int_{\tilde\Gamma_\theta}\left[ T_2\wedge \tilde\theta^*(\tilde x^0)-T_1\right] \hat\nu(d\tilde x^0)
-\int_{\tilde\Gamma_\theta}\left[ T_2\wedge \tilde\theta^*(\tilde x^0)-T_1\right]\check\eta_A([0,T_1]|\tilde x^0)\frac{a(\tilde x^0)}{n(\tilde x^0)}\hat\nu(d\tilde x^0)\\
&&-\int_{\tilde\Gamma_\theta}\int_{(T_1,T_2\wedge\tilde\theta^*(\tilde x^0)]} \left[T_2\wedge\tilde\theta^*(\tilde x^0)-u\right]\check\eta_A(du|\tilde x^0)\frac{a(\tilde x^0)}{n(\tilde x^0)}\hat\nu(d\tilde x^0)\\
&&+\int_{\bf D} w^A_{T_1,T_2,\tilde\Gamma}(\phi((\tilde x^0,0),u))\check\eta(d\tilde x^0\times du).
\end{eqnarray*}
For the last but one integral, note that $\check\eta_A(\{\tilde\theta^*(\tilde x^0)\}|\tilde x^0)=0$ for $\hat{\eta}_A$-almost all $\tilde x^0$ since $\tilde{\eta}_A$ is concentrated on $V$. The corresponding integrals over $\tilde\Gamma\setminus\tilde\Gamma_\theta$ equal zero and hence are omitted; the last term above, denoted below as ${\bf J}(\tilde{\Gamma})$, is positive. According to Lemma \ref{l4}, for $\hat\eta_A$-almost  all $\tilde x^0\in\tilde\Gamma_\theta$,
\begin{eqnarray*}
&&\int_{(T_1,T_2\wedge\tilde\theta^*(\tilde x^0)]} \left[u-T_2\wedge\tilde\theta^*(\tilde x^0)\right] \check\eta_A(du|\tilde x^0) =
\int_{(T_1,T_2\wedge\tilde\theta^*(\tilde x^0)]} \left[u-T_1\right] \check\eta_A(du|\tilde x^0)\\
&&-\left[T_2\wedge\tilde\theta^*(\tilde x^0)-T_1\right] \check\eta_A((T_1,T_2\wedge\check\theta^*(\tilde x^0)]|\tilde x^0)\\
&=&\left[T_2\wedge\tilde\theta^*(\tilde x^0)-T_1\right] \check\eta_A([T_1,T_2\wedge\tilde\theta^*(\tilde x^0)]|\tilde x^0)
-\int_{(0,T_2\wedge\tilde\theta^*(\tilde x^0)-T_1]} \check\eta_A([T_1,T_1+s)|\tilde x^0)ds\\
&&-\left[T_2\wedge\tilde\theta^*(\tilde x^0)-T_1\right] \check\eta_A((T_1,T_2\wedge\tilde\theta^*(\tilde x^0)]|\tilde x^0),
\end{eqnarray*}
so that
\begin{eqnarray}
\check\eta_\Box(Y_{T_1,T_2,\tilde\Gamma})&=& \int_{\tilde\Gamma_\theta}\left[ T_2\wedge \tilde\theta^*(\tilde x^0)-T_1\right] \hat\nu(d\tilde x^0)
-\int_{\tilde\Gamma_\theta}\left[ T_2\wedge \tilde\theta^*(\tilde x^0)-T_1\right]\check\eta_A([0,T_1]|\tilde x^0)\frac{a(\tilde x^0)}{n(\tilde x^0)}\hat\nu(d\tilde x^0)\nonumber \\
&&+\int_{\tilde\Gamma_\theta} \left[ T_2\wedge \tilde\theta^*(\tilde x^0)-T_1\right]\check\eta_A(\{T_1\}|\tilde x^0)\frac{a(\tilde x^0)}{n(\tilde x^0)}\hat\nu(d\tilde x^0)\nonumber \\
&&-\int_{\tilde\Gamma\theta}\int_{(0,T_2\wedge\tilde\theta^*(\tilde x^0)-T_1]} \check\eta_A([T_1,T_1+s)|\tilde x^0)ds~ \frac{a(\tilde x^0)}{n(\tilde x^0)}\hat\nu(d\tilde x^0)+{\bf J}(\tilde\Gamma)\nonumber \\
&=& \int_{\tilde\Gamma_\theta}\left[ T_2\wedge \tilde\theta^*(\tilde x^0)-T_1\right] \hat\nu(d\tilde x^0)
-\int_{\tilde\Gamma_\theta}\left[ T_2\wedge \tilde\theta^*(\tilde x^0)-T_1\right]\check\eta_A([0,T_1)|\tilde x^0)\frac{a(\tilde x^0)}{n(\tilde x^0)}\hat\nu(d\tilde x^0)\nonumber \\
&&-\int_{\tilde\Gamma_\theta} \int_{(T_1,T_2\wedge\tilde\theta^*(\tilde x^0)]} \check\eta_A([T_1,u)|\tilde x^0)du~\frac{a(\tilde x^0)}{n(\tilde x^0)}\hat\nu(d\tilde x^0)+{\bf J}(\tilde{\Gamma}).\label{e23}
\end{eqnarray}

According to the definitions of the measures $\tilde\eta'_\Box$ and $\check\eta'_\Box$,
\begin{eqnarray}
\check\eta'_\Box(Y_{T_1,T_2,\tilde\Gamma}) &=& \tilde\eta'_\Box(F(Y_{T_1,T_2,\tilde\Gamma}))
=\int_{\tilde V}\int_{\RR^0_+} \II\{(\tilde x^0,u)\in Y_{T_1,T_2,\tilde\Gamma}\}\left(1-\tilde p([0,u)|(\tilde{x}^0,0))\right) du~\hat\nu(d\tilde x^0)\nonumber\\
&=&\int_{\tilde{\bf V}_\nu} \int_{[T_1,T_2\wedge\tilde\theta^*(\tilde x^0))} \II\{\tilde x^0\in\tilde\Gamma\}\left(1-\tilde p([0,u)|(\tilde{x}^0,0)) \right) du~\hat\nu(d\tilde x^0)\label{e26} \\
&=& \int_{\tilde\Gamma}\int_{[T_1,T_2\wedge\tilde\theta^*(\tilde x^0))} du~\hat\nu(d\tilde x^0)
-\int_{\tilde\Gamma} \int_{[T_1,T_2\wedge\tilde\theta^*(\tilde x^0))} \left( \tilde p([0,T_1)|(\tilde{x}^0,0))\nonumber\right.\\
&&\left.+\tilde p([T_1,u)|(\tilde{x}^0,0))\right) du~ \hat\nu(d\tilde x^0)\nonumber\\
&=& \int_{\tilde\Gamma_\theta}\left[ T_2\wedge\tilde\theta^*(\tilde x^0)-T_1\right] \hat\nu(d\tilde x^0)
-\int_{\tilde\Gamma_\theta} \tilde p([0,T_1)|(\tilde{x}^0,0))\left[ T_2\wedge\tilde\theta^*(\tilde x^0)-T_1\right] \hat\nu(d\tilde x^0)\nonumber \\
&&-\int_{\tilde\Gamma_\theta} \int_{(T_1,T_2\wedge\tilde\theta^*(\tilde x^0)]} \tilde p([T_1,u)|(\tilde{x}^0,0)) du~\hat\nu(d\tilde x^0). \nonumber
\end{eqnarray}

Since $Y_{T_1,T_2,\tilde\Gamma\setminus\tilde\Gamma_\theta}=\emptyset$,
\begin{eqnarray*}
\check\eta_\Box(Y_{T_1,T_2,\tilde\Gamma\setminus\tilde\Gamma_\theta})-\check\eta'_\Box(Y_{T_1,T_2,\tilde\Gamma\setminus\tilde\Gamma_\theta})=0.
\end{eqnarray*}

It remains to consider the set $\tilde\Gamma_\theta$.
Below, we split it  into three measurable subsets:
\begin{eqnarray*}
\tilde\Gamma_1 &\defi & \tilde\Gamma_\theta\cap\{\tilde x^0: ~u^*(\tilde x^0)< T_1\},\\
\tilde\Gamma_2 &\defi & \tilde\Gamma_\theta\cap\{\tilde x^0: ~u^*(\tilde x^0)\ge T_2\wedge\tilde\theta^*(\tilde x^0)\},\\
\mbox{ and}~~\tilde\Gamma_3 &\defi & \tilde\Gamma_\theta\cap\{\tilde x^0: ~T_1\le u^*(\tilde x^0)<T_2\wedge\tilde\theta^*(\tilde x^0)\}.
\end{eqnarray*}

For  each $\tilde x^0\in\tilde\Gamma_1$, $\tilde p([0,u)|(\tilde{x}^0,0))=1$ for all $u\in[T_1,T_2\wedge\tilde\theta^*(\tilde x^0))$.  Hence, according to (\ref{e26}) with $\tilde{\Gamma}=\tilde{\Gamma}_1$, $\check\eta'_\Box(Y_{T_1,T_2,\tilde\Gamma_1})=0$ and
\begin{eqnarray*}
\check\eta_\Box(Y_{T_1,T_2,\tilde\Gamma_1})-\check\eta'_\Box(Y_{T_1,T_2,\tilde\Gamma_1})\ge 0.
\end{eqnarray*}

For each $\tilde x^0\in\tilde\Gamma_2$ (see the point $\tilde x^0_1$ on Figure  \ref{fig5}),
$$\tilde p([0,T_1)|(\tilde{x}^0,0))=\check\eta_A([0,T_1)|\tilde x^0)\frac{a(\tilde x^0)}{n(\tilde x^0)}~\mbox{ and }
\tilde p([T_1,u)|(\tilde{x}^0,0))=\check\eta_A([T_1,u)|\tilde x^0)\frac{a(\tilde x^0)}{n(\tilde x^0)}$$
for all $u\in(T_1,T_2\wedge\tilde\theta^*(\tilde x^0)]$. Therefore, by (\ref{e23}) and (\ref{e26}),
$$\check\eta_\Box(Y_{T_1,T_2,\tilde\Gamma_2})-\check\eta'_\Box(Y_{T_1,T_2,\tilde\Gamma_2})={\bf J}(\tilde{\Gamma}_2)\ge 0.$$

For the set $\tilde\Gamma_3$ (the typical points in $\tilde\Gamma_3$ are $\tilde x^0_2$ and $\tilde x^0_3$ on Figure \ref{fig5}), we compute $\check\eta_\Box(Y_{T_1,T_2,\tilde\Gamma_3})$ and $\check\eta'_\Box(Y_{T_1,T_2,\tilde\Gamma_3})$ using the representation $Y_{T_1,T_2,\tilde\Gamma_3}=Y^1\cup Y^2$, where
$$Y^1\defi\{(\tilde x^0, u):~\tilde x^0\in\tilde\Gamma_3,~T_1\le u<u^*(\tilde x^0)\};~~~Y^2\defi\{(\tilde x^0,u):~\tilde x^0\in\tilde\Gamma_3,~u^*(\tilde x^0)\le u< T_2\wedge \tilde\theta^*(\tilde x^0)\}.$$
To compute $\check\eta_\Box(Y^1)$, we introduce the function
\begin{eqnarray*}
&&w(y)=w((\tilde y,t))
= \II\{h(y)\in\tilde\Gamma_3\}\times \left\{\begin{array}{ll}
u^*(h(y))-T_1, & \mbox{ if } t\le T_1;\\
u^*(h(y))-t, & \mbox{ if } T_1<t\le u^*(h(y));\\
0& \mbox { if } t>u^*(h(y))
\end{array}\right.
\end{eqnarray*}
(cf (\ref{e27})). Calculations similar to those presented above, lead to the following version of expression (\ref{e23}):
\begin{eqnarray*}
\check\eta_\Box(Y^1)&=& \int_{\tilde\Gamma_3}\left[ u^*(\tilde x^0)-T_1\right] \hat\nu(d\tilde x^0)
-\int_{\tilde\Gamma_3}\left[ u^*(\tilde x^0)-T_1\right]\check\eta_A([0,T_1)|\tilde x^0)\frac{a(\tilde x^0)}{n(\tilde x^0)}\hat\nu(d\tilde x^0) \\
&&-\int_{\tilde\Gamma3} \int_{(T_1,u^*(\tilde x^0)]} \check\eta_A([T_1,u)|\tilde x^0)du~\frac{a(\tilde x^0)}{n(\tilde x^0)}\hat\nu(d\tilde x^0)+{\bf J}^1.
\end{eqnarray*}
The last term is similar to $\bf J(\tilde{\Gamma})$, its calculation is based on the function similar to $w^A_{T_1,T_2,\tilde\Gamma}$: one only has to replace $\tilde\Gamma$ with $\tilde\Gamma_3$ and $\tilde\theta^*(\cdot)$ with $u^*(\cdot)$. Like previously, ${\bf J}^1\ge 0$. Again, similarly to (\ref{e26}), we have
\begin{eqnarray*}
\check\eta'_\Box(Y^1) &=&
\int_{\tilde\Gamma_3}\left[u^*(\tilde x^0)-T_1\right] \hat\nu(d\tilde x^0)
-\int_{\tilde\Gamma_3} \left[u^*(\tilde x^0)-T_1\right]\tilde p([0,T_1)|(\tilde x^0,0)) \hat\nu(d\tilde x^0) \\
&&-\int_{\tilde\Gamma_3} \int_{(T_1,u^*(\tilde x^0)]} \tilde p([T_1,u)|(\tilde x^0,0)) du~\hat\nu(d\tilde x^0)
\end{eqnarray*}
and, like in the case of $\tilde\Gamma_2$, for each $\tilde{x}_0\in\tilde{\Gamma}_3$
$$\tilde p([0,T_1)|(\tilde x^0,0))=\check\eta_A([0,T_1)|\tilde x^0)\frac{a(\tilde x^0)}{n(\tilde x^0)}~\mbox{ and }
\tilde p([T_1,u)|(\tilde x^0,0))=\check\eta_A([T_1,u)|\tilde x^0)\frac{a(\tilde x^0)}{n(\tilde x^0)}$$
for all $u\in(T_1,u^*(\tilde x^0)]$. Therefore,
$$\check\eta_\Box(Y^1)-\check\eta'_\Box(Y^1)={\bf J}^1\ge 0.$$
Finally, similarly to (\ref{e26}),
$$\check\eta'_\Box(Y^2)=\int_{\tilde\Gamma_3}\int_{[u^*(\tilde x^0),T_2\wedge\tilde\theta^*(\tilde x^0))} (1-\tilde p([0,u)|(\tilde x^0,0)))du~\hat\nu(d\tilde x^0)=0$$
because for each $\tilde{x}^0\in\tilde{\Gamma}_3$, $\tilde p([0,u)|(\tilde x^0,0))=1$ for all $u>u^*(\tilde x^0)$. Hence,
$$\check\eta_\Box(Y^2)-\check\eta'_\Box(Y^2)\ge 0.$$

To summarize, $\check\eta_\Box(Y_{T_1,T_2,\tilde\Gamma_3})-\check\eta'_\Box(Y_{T_1,T_2,\tilde\Gamma_3})\ge 0$, and thus $\check\eta_\Box(Y_{T_1,T_2,\tilde\Gamma})-\check\eta'_\Box(Y_{T_1,T_2,\tilde\Gamma})\ge 0$ for all $\tilde{\Gamma}\in {\cal B}(\tilde{V})$ and $0\le T_1<T_2<\infty.$

Therefore,  $\check\eta_\Box\ge\check\eta'_\Box$ set-wise on $\bf D$, and hence $\tilde\eta_\Box\ge \tilde\eta'_\Box$ on $V$. Since the  measures $\tilde\eta_\Box$ and $\tilde\eta'_\Box$ are both normal,  the difference $\tilde\eta_\Box-\tilde\eta'_\Box$ is a normal  measure on $V$ by Lemma \ref{l101}(b).

The proof is completed. \hfill$\Box$
\bigskip

\par\noindent\underline{Proof of Theorem \ref{t2}.}
When $x\in V^c$, we fix $\pi_i(d\theta\times da)\defi \delta_{f^*(x)}(d\theta\times da)$, where  $f^*(x)=(\infty,\hat a)$ as usual.  Below, for two finite or normal measures $\zeta^1$ and $\zeta^2$ on $V$, the inequality $\zeta^1(dx)\le \zeta^2(dx)$ is understood set-wise. The same concerns measures on $V\times{\bf A}$.

Let $p'_A(da|x)$ be the stochastic kernel on $\bf A$ given $V$ coming from the decomposition $\eta(dx\times da)=p'_A(da|x)\eta(dx\times{\bf A})$. For all $i\ge 1$, we put
\begin{eqnarray*}
p^i_A(da|x,\theta)\equiv p'_A(da|\phi(x,\theta))
\end{eqnarray*}
for  $x\in V$, $\theta<\theta^*(x)$, and $p^i_A(da|x,\theta)$ is an arbitrarily fixed stochastic kernel on $\bf A$ for $x\in V$, $\theta\ge\theta^*(x)$.

We will prove by induction the following statement.

For each $i\ge 1,$ there is a stochastic kernel $\pi_i$ on ${\bf B}=\bar\RR^0_+\times{\bf A}$ given $V$, having the form
\begin{eqnarray*}
\pi_i(d\theta\times da|x)=p^i_T(d\theta|x)p^i_A(da|x,\theta),
\end{eqnarray*}
such that, for each $n\ge 1$ and the sequence $\{\pi_i\}_{i=1}^n$, the following assertions are fulfilled.

(i) $p_T^i([\theta^\ast(x),\infty)|x)=0$ for $x\in V,$ $i=1,2,\dots,n,$ and the (partial) aggregated occupation measures $\{\tilde\eta^i\}_{i=0}^n$, defined as in Lemma \ref{l11}, exhibit the following properties:
\begin{eqnarray*}
\tilde\eta^n(dx\times\Box)&\le & \eta(dx\times\Box)~~~~~\mbox{ and }\\
\tilde\eta^n(dx\times da)= \tilde\eta^n(dx\times{\bf A}) p'_A(da|x) &\le & \eta(dx\times{\bf A})p'_A(da|x)=\eta(dx\times da) ~\mbox{on ${\cal B}(V\times\textbf{A})$}.
\end{eqnarray*}

(ii) The measure $\nu^n(dx):=P^{\pi}_{x_0}(X_n\in dx)$ on $V$ is such that, for each function $w\in{\bf W}$,
\begin{eqnarray}
0&=& \int_{V} w(x)\nu^n(dx)+\int_{V} \chi w(x)[\eta-\tilde\eta^n](dx\times\Box)-\int_{V} w(x)[\eta-\tilde\eta^n](dx\times{\bf A}) \nonumber\\
&&+\int_{V\times{\bf A}} w(l(x,a))[\eta-\tilde\eta^n](dx\times da), \label{en8}
\end{eqnarray}
and all the integrals here are finite. Note that $\nu^n$ is uniquely defined by the finite sequence $\{\pi_i\}_{i=1}^n$: see (\ref{e33}); moreover, $\nu^n(\tilde V\times\{t:~t>0\})=0$.

After that, $\pi^\eta:=\{\pi_i\}_{i=1}^\infty$ will be the desired Markov strategy.

When $n=0$, $\tilde\eta^0(dy\times\Box)\equiv 0$,  $\tilde\eta^0(dy\times da)\equiv 0$, and $\nu^0(dx)=\delta_{x_0}(dx)$. Assertions (i) and (ii) are obviously fulfilled because the normal measure $\eta$ satisfies equation (\ref{e17}).

Suppose assertions (i) and (ii) hold true for $i=0,1,2,\ldots,n\ge 0$.  We apply Lemma \ref{l8} to the measures $\nu\defi \nu^n$, $\tilde\eta_\Box\defi (\eta-\tilde\eta^n)(dx\times\Box)$, $\tilde\eta_A\defi (\eta-\tilde\eta^n)(dx\times{\bf A})$, and $\tilde\eta\defi(\eta-\tilde\eta^n)(dx\times da)$ satisfying equation (\ref{en8}). All of them are finite, maybe apart from $\tilde\eta_\Box$, which is normal by Lemma \ref{l101}(b) and Lemma \ref{l11}.
As a result, we have the stochastic kernel $\tilde p(dt|x)$ on $\bar\RR^0_+$ given $V$ and the measures
\begin{equation}\label{en2}
\tilde\eta'_A(dx)\le (\eta-\tilde\eta^n)(dx\times{\bf A})~~\mbox{ and }~~ \tilde\eta'_\Box(dx)\le (\eta-\tilde\eta^n)(dx\times\Box)
\end{equation}
on $V$, which satisfy equation (\ref{e13}):
\begin{equation}\label{en4}
0=\int_{V} w(x)\nu^n(dx)+\int_{V} \chi w(x)\tilde\eta'_\Box(dx)-\int_{V} w(x)\tilde\eta'_A(dx),~~~w\in{\bf W}.
\end{equation}
All the integrals here are finite.

For $x\in V$, we put
\begin{eqnarray*}
p^{n+1}_T(d\theta|x)\defi \tilde p(d\theta|x).
\end{eqnarray*}
Then by Lemma \ref{l8}, $p^{n+1}_T([\theta^\ast(x),\infty)|x)=0$ for all $x\in V.$
All the kernels $\{\pi_i\}_{i=1}^n$ were built on the previous steps of the induction. According to the definition of the measure $\tilde\eta^{n+1}$,
\begin{eqnarray}
\tilde\eta^{n+1}(\Gamma\times\Box) &=& \tilde\eta^n(\Gamma\times\Box)+\int_{V}\int_{\RR^0_+} \delta_{\phi(x,u)}(\Gamma) \tilde p([u,\infty]|x)du~\nu^n(dx) \nonumber\\
&=&\tilde\eta^n(\Gamma\times\Box)+\tilde\eta'_\Box(\Gamma)\le\eta(\Gamma\times\Box), ~~~~~~~~~~\Gamma\in{\cal B}(V);\label{en6} \\
\tilde\eta^{n+1}(\Gamma\times{\bf A}) &=& \tilde\eta^n(\Gamma\times{\bf A})+\int_{V}\int_{\RR^0_+} \delta_{\phi(x,u)}(\Gamma) \tilde p(du|x)\nu^n(dx)=\tilde\eta^n(\Gamma\times{\bf A})+\tilde\eta'_A(\Gamma)\le\eta(\Gamma\times{\bf A}),\nonumber \\
&&~~~~~~~~\Gamma\in{\cal B}(V). \label{en7}
\end{eqnarray}
Inequalities are valid according to the basic properties of the measures $\tilde\eta'_\Box$ and $\tilde\eta'_A$ presented in (\ref{en2}).
Recall that
\begin{eqnarray*}
\tilde\eta^{n+1}(\Gamma_X\times\Gamma_A) &=& \tilde\eta^n(\Gamma_X\times\Gamma_A)+\int_{V}\int_{\RR^0_+} \delta_{\phi(x,\theta)}(\Gamma_X)
p^{n+1}_A(\Gamma_A|x,\theta) p^{n+1}_T(d\theta|x)\nu^n(dx),\\
&&~~~~~~\Gamma_X\in{\cal B}(V),~\Gamma_A\in{\cal B}({\bf A}).
\end{eqnarray*}
Since $\nu^n(\tilde V\times\{t:~t>0\})=0$, the last term equals
\begin{equation}\label{en10}
I\defi\int_{\tilde V}\int_{[0,\tilde\theta^*(\tilde x^0))} \delta_{\phi((\tilde x^0,0),\theta)}(\Gamma_X) p'_A(\Gamma_A|\phi((\tilde x^0,0),\theta)) \tilde p(d\theta|(\tilde x^0,0)) {\hat\nu}^n(d\tilde x^0),
\end{equation}
where $\hat{\nu}^n(\Gamma)\defi\nu^n(\{(\tilde x^0,0),~\tilde x^0\in\Gamma\})$. According to Lemma \ref{l8}, for all $\Gamma\in{\cal B}({\bf D})$ and for the mapping $F$ as in Definition \ref{d7},
\begin{eqnarray*}
\tilde\eta'_A(F(\Gamma))&=& \int_{\tilde V}\int_{\RR^0_+} \II\left\{\phi((\tilde x^0,0),u)\in\{y=\phi((\tilde x^0,0),t):~(\tilde x^0,t)\in\Gamma\}\right\} \tilde p(du|(\tilde x^0,0)) \hat\nu^n(d\tilde x^0)\\
&=& \int_{\tilde V}\int_{\RR^0_+} \II\{(\tilde x^0,u)\in\Gamma\}\tilde p(du|(\tilde x^0,0))\hat\nu^n(d\tilde x^0)=\check{\tilde\eta}'_A(\Gamma).
\end{eqnarray*}
Lemma \ref{l9}(a) implies that, for each bounded measurable function $g$ on $V$,
\begin{equation}\label{en5}
\int_{V} g(x)\tilde\eta'_A(dx)=\int_{\bf D}g(\phi((\tilde x^0,0),u)) \check{\tilde\eta}'_A(d\tilde x^0\times du)
=\int_{\tilde V}\int_{[0,\tilde\theta^*(\tilde x^0))} g(\phi((\tilde x^0,0),u)) \tilde p (du|(\tilde x^0,0))\hat\nu^n(d\tilde x^0).
\end{equation}
Therefore, for each $\Gamma_X\in{\cal B}(V)$,
$$I=\int_{V} \delta_x(\Gamma_X)p'_A(\Gamma_A|x)\tilde\eta'_A(dx)=\int_{\Gamma_X} p'_A(\Gamma_A|x)\tilde\eta'_A(dx),$$
meaning that on ${\cal B}(V\times \textbf{A})$
\begin{eqnarray*}
\tilde\eta^{n+1}(dx\times da)&=& \tilde\eta^n(dx\times da)+p'_A(da|x)\tilde \eta'_A(dx)=\tilde\eta^n(dx\times{\bf A})p'_A(da|x)+\tilde \eta'_A(dx)p'_A(da|x)\\
&=& \tilde\eta^{n+1}(dx\times{\bf A}) p'_A(da|x)\\
&\le &\tilde\eta^n(dx\times{\bf A}) p'_A(da|x)+[\eta-\tilde\eta^n](dx\times{\bf A}) p'_A(da|x)=\eta(dx\times{\bf A}) p'_A(da|x).
\end{eqnarray*}
The second equality is by the inductions supposition, the third equality follows from (\ref{en7}), and the inequality is according to the basic property (\ref{en2}) of the measure $\tilde\eta'_A$.

Property (i) for $n+1$ is established, recall also inequality (\ref{en6}).

For the proof of Item (ii), note that, by (\ref{en8}) at $n$, (\ref{en4}), (\ref{en6}), and (\ref{en7}), we have equation
\begin{eqnarray*}
0&=&\int_{V} \chi w(x)[\eta-\tilde\eta^{n+1}](dx\times\Box)-\int_{V} w(x)[\eta-\tilde\eta^{n+1}](dx\times{\bf A})\\
&&+\int_{V\times{\bf A}} w(l(x,a))[\eta-\tilde\eta^{n+1}](dx\times da)+\int_{V\times{\bf A}} w(l(x,a))[\tilde\eta^{n+1}-\tilde\eta^{n}](dx\times da)
\end{eqnarray*}
valid for all functions $w\in{\bf W}$, and all the integrals here are finite. According to property (i) for $n$ and $n+1$, the stochastic kernel $p'_A(da|x)$ is the same in the decompositions $\tilde\eta^n(dx\times da)=\tilde\eta^n(dx\times{\bf A})p'_A(da|x)$ and  $\tilde\eta^{n+1}(dx\times da)=\tilde\eta^{n+1}(dx\times{\bf A})p'_A(da|x)$. Thus, the last integral, according to (\ref{en7}), equals
$$\int_{V}\int_{\bf A} w(l(x,a)) p'_A(da|x)\tilde\eta'_A(dx),$$
i.e., the function $w$ is integrated with respect to the measure
$$m(\Gamma)=\int_{V}\int_{\bf A} \delta_{l(x,a)}(\Gamma) p'_A(da|x)\tilde\eta'_A(dx),~~~\Gamma\in{\cal B}(V),$$
and it remains to show that this measure coincides with $\nu^{n+1}$ on $V$.

From equation (\ref{en5}), we have for all $\Gamma\in{\cal B}(V)$:
$$m(\Gamma)=\int_{\tilde V}\int_{[0,\tilde\theta^*(\tilde x^0))}\int_{\bf A} \delta_{l(\phi((\tilde x^0,0),u),a)} (\Gamma) p'_A(da|\phi((\tilde x^0,0),u)) p^{n+1}_T(du|(\tilde x^0,0))\hat{\nu}^n(d\tilde x^0),$$
and, keeping in mind that $\nu^n(\tilde V\times\{t:~t>0\})=0$, we have from (\ref{e33}):
$$\nu^{n+1}(\Gamma)=\int_{\tilde V}\int_{[0,\tilde\theta^*(\tilde x^0))}\int_{\bf A} \delta_{l(\phi((\tilde x^0,0),\theta),a)} (\Gamma) p^{n+1}_A(da|(\tilde x^0,0),\theta) p^{n+1}_T(d\theta|(\tilde x^0,0)) \hat\nu^n(d\tilde x^0)=m(\Gamma)$$
for all $\Gamma\in{\cal B}(V)$ because $p^{n+1}_A=p'_A$.

The proof of the induction statement for $n+1$ is completed.

According to Lemma \ref{l11}, for the constructed Markov strategy $\pi=\{\pi_i\}_{i=1}^\infty$ and for the corresponding aggregated occupation measure $\tilde\eta$, we have the convergence $\tilde\eta^n\uparrow\tilde\eta$ set-wise as $n\to\infty$. Since $\tilde\eta^n\le\eta$ set-wise on $V\times{\bf A}_\Box$, the desired set-wise inequality $\tilde\eta\le\eta$ follows.

All the properties enlisted in Definition \ref{JulyRem01} are obviously satisfied for the strategy $\pi$.  \hfill$\Box$
\bigskip

\par\noindent\underline{Proof of Corollary \ref{corol1}.}  We denote by $Val(\ref{e106})$ and $Val(\ref{e123})$ the minimal values of linear programs (\ref{e106})  and (\ref{e123}), respectively. Recall that linear program (\ref{e106}) has an optimal solution by Proposition \ref{pr1}.

Suppose the finite measure $\mu^*$ on $V\times\bar\RR^0_+\times{\bf A}$ (concentrated on $\textbf{M}\times\textbf{A}$) solves linear program (\ref{e106}). Then the aggregated occupation measure $\eta^*$, induced by $\mu^*$, is normal by Lemma \ref{l101}(a) and  satisfies equation (\ref{e17}) according to Theorem \ref{t1}. The constraints-inequalities  are also fulfilled by $\eta^*$. Thus
\begin{eqnarray*}
\infty>Val(\ref{e106})=\int_{V\times{\bf A}_\Box} C_0(x,a)\eta^*(dx\times da)\ge Val(\ref{e123}).
\end{eqnarray*}
In case the last inequality is strict, there exists a feasible solution $\eta$ to linear program (\ref{e123}) satisfying inequality
\begin{eqnarray*}
\int_{V\times{\bf A}_\Box} C_0(x,a)\eta(dx\times da)<Val(\ref{e106}).
\end{eqnarray*}
Consider the induced reasonable Markov strategy $\pi^\eta$ as in
Theorem \ref{t2} and the corresponding aggregated occupation measure $\tilde\eta$.  Since $C_j\ge 0$ for $j=0,1,\ldots, J$, all the conditions in linear program (\ref{e106}) are satisfied for $\mu^{\pi^\eta}$ and
\begin{eqnarray*}
\int_{V\times{\bf A}_\Box} C_0(x,a)\tilde\eta(dx\times da)\le\int_{V\times{\bf A}_\Box} C_0(x,a)\eta(dx\times da)<Val(\ref{e106})<\infty.
\end{eqnarray*}
The measure $\mu^{\pi^\eta}$ cannot take infinite value as explained above linear program (\ref{SashaLp02}). We obtained a contradiction to the optimality of the measure $\mu^*$. Hence, $Val(\ref{e106})=Val(\ref{e123})$, and the measure $\eta^*$ solves linear program (\ref{e123}).

Suppose now that the measure $\eta^*$ on $V\times{\bf A}_\Box$ solves linear program (\ref{e123}) and consider the reasonable Markov strategy $\pi^*=\pi^{\eta^*}$ as in Theorem \ref{t2}. The corresponding occupation measure $\mu^{\pi^*}$ is feasible in linear program (\ref{e106}). More detailed reasoning is similar to that presented above. Therefore, for the aggregated occupation measure $\tilde\eta$ induced by $\mu^{\pi^*}$, we have relations
\begin{eqnarray*}
Val(\ref{e106})\le \int_{V\times{\bf A}_\Box} C_0(x,a)\tilde\eta(dx\times da)\le  \int_{V\times{\bf A}_\Box} C_0(x,a)\eta^*(dx\times da)=Val(\ref{e123}).
\end{eqnarray*}
But we have shown that $Val(\ref{e123})=Val(\ref{e106})$,  so that
\begin{eqnarray*}
\int_{V\times{\bf A}_\Box} C_0(x,a)\tilde\eta(dx\times da)=Val(\ref{e106})
\end{eqnarray*}
meaning that the measure $\mu^{\pi^*}$ solves linear program (\ref{e106}).
The proof is completed. \hfill$\Box$

\section{Acknowledgement}

This research was supported by the Royal Society International Exchanges award IE160503. We would like to  thank Prof.A.Plakhov for his initial participation in this work and for his proof of Lemma \ref{l1}.

\appendix
\section{Appendix} \label{s5}

Lemma \ref{l1} and its proof presented below are similar to Lemma 2.2 in \cite{b2}, where the authors assumed that $E$ was a subset of an  Euclidean space.

Let $E$ be an arbitrary set and $\phi : E \times \RR_+^0 \to E$ be a flow in $E$ possessing the semigroup property.

\begin{definition}\label{JulyDef01}
A function $w : E \to \RR$ is said to be absolutely continuous along the flow if for all $x \in E$ the function $t \mapsto w(\phi(x,t)), \ t \in \RR_+^0$ is absolutely continuous. It is called increasing (decreasing) along the flow if so is the function $t\to w(\phi(x,t))$, $t\in\RR^0_+$ for all $x\in E$.
\end{definition}

\begin{lemma}\label{l1} Suppose function
$w$ is absolutely continuous along the flow $\phi$. Then the following assertions are valid.

(a) There exists a function $\chi w : E \to \RR$ such that, for any $x\in E$, the function $\chi w(\phi(x,s))$ is Lebesgue integrable with respect to $s$ on any finite interval $[0,t]\subset \RR_+^0$ and
\begin{equation}\label{barrow}
w(\phi(x,t)) - w(x) = \int_{[0,t]} \chi w(\phi(x,s))\, ds
\end{equation}
for all $x \in E$ and $t \ge 0$.

(b) If, additionally, $E$ is a measurable space (that is, is equipped with a $\sigma$-algebra of subsets), $w$ is measurable, and the functions $\phi(\cdot,t) : E \to E$ are measurable for all $t \ge 0$, then the function $\chi w$ satisfying (a) can be chosen measurable.
\end{lemma}

\par\noindent\underline{Proof.} We provide one common proof for (a) and (b) underlining the measurability properties as soon as they appear.

Define the functions
$$
\overline{W}(x) := \overline{\lim}_{n\to\infty} \frac{w(\phi(x, \frac{1}{n})) - w(x)}{1/n}, \qquad
\underline{W}(x) := \underline{\lim}_{n\to\infty} \frac{w(\phi(x, \frac{1}{n})) - w(x)}{1/n}
$$
and the set $D := \{ x \in E : \ \overline{W}(x) = \underline{W}(x) \ne \pm\infty \}$. Let us additionally define the function $W : D \to \RR$ by $W(x) := \overline{W}(x)$; that is, $W(x)$ coincides with the limit $\lim_{n\to\infty} n\, [w(\phi(x, \frac{1}{n})) - w(x)]$, if it exists and is finite.

If $w$ and $\phi(\cdot,t)$ are measurable, then $w(\phi(x, \frac{1}{n}))$ is also measurable. Hence the functions $\overline{W}$ and $\underline{W}$ are measurable as the upper and lower limits of the sequence of measurable functions $n\,[w(\phi(x, 1/n)) - w(x)]$. Consequently, the set $D$ is also measurable.

Define the function $\chi w$ on $E$ by
\begin{equation}\label{ea1}
\chi w(x) :=
\left\{
\begin{array}{ll}
W(x), & \text{if} \ \ x \in D;\\
g(x), & \text{otherwise},
\end{array}
\right.
\end{equation}
where $g$ is any function. In the measurable case we take $g$ to be measurable and readily get that $\chi w$ is also measurable.

Since $w$ is absolutely continuous along the flow then for any $x \in E$ there exists a subset of full measure $T_x \subset \RR_+ $ such that the derivative $\frac{d}{dt}\, w(\phi(x,t))$ exists and is finite for all values $t \in T_x$. For any such value (let it now be denoted by $s \in T_x$) we can write down the following (below we denote $x' = \phi(x,s)$ and use the semigroup property of the flow)
$$
\frac{dw(\phi(x,t))}{dt}\Big\rfloor_{t = s} = \lim_{\varepsilon\to 0} \frac{w(\phi(x, s + \varepsilon)) - w(\phi(x, s))}{\varepsilon}
= \lim_{n\to\infty} \frac{w(\phi(x', \frac{1}{n})) - w(x')}{1/n}.
$$
The latter value exists and is finite, and therefore coincides with $W(x')$. This argument also shows that $\phi(x, T_x) \subset D$.

Since $w$ is absolutely continuous along the flow, one can write down
$$
w(\phi(x,t)) - w(x) = \int_{[0,t]\cap T_x} \frac{dw(\phi(x,\tau))}{d\tau}\Big\rfloor_{\tau = s}\, \, ds = \int_{[0,t]\cap T_x} W(\phi(x,s))\, ds.
$$
Now taking into account that $[0,\, t] \setminus T_x$ has Lebesgue measure zero and $\chi w$ is an extension of $W$ to $E$, we conclude that the latter integral coincides with $\int_{[0,t]} \chi w(\phi(x,s))\, ds$, and so, formula (\ref{barrow}) is proved. \hfill $\Box$
\bigskip

\par\noindent\underline{Proof of Lemma \ref{la5}.} In this proof, let us denote by $\rho$ and $\tilde{\rho}$ the compatible metrics on $\tilde{\bf X}\times\RR^0_+$ and $\tilde{\bf X}.$ If $y_n\to y$, where $y_n=(\tilde y_n,t_n)$, $y=(\tilde y,t)\in{\bf X}$, then the sequence $\{t_n\}_{n=1}^\infty$ is bounded: $t,t_n\in[0,T]$ for some $T<\infty$. Now
$\tilde\rho(h(y_n),h(y))\le \sup_{t\in[0,T]} d(t)\rho(y_n,y)\to 0$. Thus, $h$ is continuous. The continuity of the mapping $h$ and of the original flow $\tilde\phi$ immediately implies that the flows $\tilde\phi$ and $\phi$ in the reverse time are continuous.

The mapping $F$ is continuous because the flow $\tilde\phi$ is continuous. It is a bijection from $\tilde{\bf X}\times \RR^0_+$ to ${\bf X}$, and the inverse mapping $F^{-1}(y)=(h(y),\tau_y)$ is continuous, as has been proved above. (For $y=(\tilde y,t)\in{\bf X}$, $\tau_y=t$ is obviously a continuous function of $y$.) Thus, $F$ is a homeomorphism, and $\bf X$ is a Borel space, being the homeomorphic image of the Borel space $\tilde{\bf X}\times \RR^0_+$. See also \cite[Prop.7.15]{Bertsekas:1978}.
\hfill$\Box$
\bigskip

\par\noindent\underline{Proof of Lemma \ref{l101}.} (a) The measure $\eta$ is finite on $V\times{\bf A}$ because the measure $\mu$ is finite. Recall that the measure $\mu$ is concentrated on $\tilde V\times\{0\}\times\bar\RR^0_+\times{\bf A}$. For the measure $\eta(dx\times\Box)$ on $V$, we have
\begin{eqnarray*}
\check\eta(\Gamma)&=&\eta(F(\Gamma)\times\Box)=\int_{\RR^0_+}\left\{\int_{\tilde V} \II\{\phi((\tilde x,0),u)\in F(\Gamma)\} \mu(d\tilde x\times\{0\}\times[u,\infty]\times{\bf A})\right\}du\\
&=&\int_{\RR^0_+}\left\{\int_{\tilde V} \II\{(\tilde x,u)\in \Gamma\} \mu(d\tilde x\times\{0\}\times[u,\infty]\times{\bf A})\right\}du\\
&=& \int_{\tilde V}\int_{\RR^0_+}\II\{(\tilde x,u)\in\Gamma\}\int_{\bf A} p_T([u,\infty]|(\tilde x,0),a) p_A(da|(\tilde x,0))du~\mu(d\tilde x\times\{0\}\times\bar\RR^0_+\times{\bf A}),
\end{eqnarray*}
for all $\Gamma\in{\cal B}({\bf D})$: see (\ref{e14p}) and (\ref{e10}). Thus, $\eta$ on $V\times \textbf{A}_\Box$ is normal.

(b) As mentioned in the proof of part (a), $\eta(V\times{\bf A})<\infty$.

Consider the measures $\zeta^1(dy):=\eta^1(dy\times\Box)$,  $\zeta^2(dy):=\eta^2(dy\times\Box)$,  and
$\zeta(dy):=\eta(dy\times\Box)$. Since $\zeta\ge 0$, we have $\check\zeta=\check\zeta^1-\check\zeta^2\ge 0$ as well.
If
$$\check\zeta^1(d\tilde x^0\times du)=g^1(\tilde x^0,u)du~L^1(d\tilde x^0)~\mbox{ and } \check\zeta^2(d\tilde x^0\times du)=g^2(\tilde x^0,u)du~L^2(d\tilde x^0),$$
then we put $L:=L^1+L^2$ and
$$g(\tilde x^0,u):=\left[ g^1(\tilde x^0,u)\frac{dL^1}{dL}(\tilde x^0)-g^2(\tilde x^0,u)\frac{dL^2}{dL}(\tilde x^0)\right]
\II\left\{g^1(\tilde x^0,u)\frac{dL^1}{dL}(\tilde x^0)-g^2(\tilde x^0,u)\frac{dL^2}{dL}(\tilde x^0)\ge 0\right\}.$$
The measurable set
$$\Gamma:=\left\{(\tilde x^0,u):~\tilde x^0\in\tilde V,~0\le u<\tilde\theta^*(\tilde x^0),~g^1(\tilde x^0,u)\frac{dL^1}{dL}(\tilde x^0)-g^2(\tilde x^0,u)\frac{dL^2}{dL}(\tilde x^0)< 0\right\}\subset{\bf D}$$
is null with respect to the measure $L(d\tilde x^0)\times du$ because, otherwise, we would have for some $t<\infty$, $\int_{\Gamma_t} L(d\tilde{x}^0)\times du>0$ for the set
\begin{eqnarray*}
\Gamma_t:=\Gamma\cap\{(\tilde x^0,u):~\tilde x^0\in\tilde V,~u\le t\},
\end{eqnarray*}
and yield a desired contradiction:
\begin{eqnarray*}
0>\int_{\Gamma_t} \left[g^1(\tilde x^0,u)\frac{dL^1}{dL}(\tilde x^0)-g^2(\tilde x^0,u)\frac{dL^2}{dL}(\tilde x^0)\right] L(d\tilde x^0)\times du=\check\zeta^1(\Gamma_t)-\check\zeta^2(\Gamma_t)=\check\zeta(\Gamma_t)\ge 0 \end{eqnarray*}
with all the terms being finite.
Now, for $\zeta=\zeta^1-\zeta^2$, we have
$$\check\zeta(d\tilde x^0\times du)=\check\zeta^1(d\tilde x^0\times du)-\check\zeta^2(d\tilde x^0\times du)=g(\tilde x^0,u)du~L(d\tilde x^0),$$
and the proof is completed.
\hfill $\Box$

\begin{lemma}\label{l4}
Suppose $m$ is a finite measure on $\RR$. Then, for each $\tau,t\in\RR$,
$$tm([\tau,\tau+t])=\int_{(0,t]} m([\tau,\tau+s))ds+\int_{(\tau,\tau+t]} (s-\tau)~dm(s)$$
and
$$tm([\tau,\tau+t))=\int_{[0,t)} m([\tau,\tau+s))ds+\int_{[\tau,\tau+t)} (s-\tau)~dm(s).$$
\end{lemma}

\par\noindent\underline{Proof.} For all cadlag (i.e., right-continuous with left limits) real-valued functions $U$ and $V$ on $\RR$ with finite variation (on finite intervals),
\begin{equation}\label{estarr}
U(t_2)V(t_2)=U(t_1)V(t_1)+\int_{(t_1,t_2]} U(s-)~dV(s)+\int_{(t_1,t_2]} V(s)~dU(s)
\end{equation}
for any $-\infty< t_1<t_2<\infty$. (See \cite[Appendix A4,\S 2]{b6}.) Equivalently, in the symmetric form:
$$U(t_2)V(t_2)=U(t_1)V(t_1)+\int_{(t_1,t_2]} U(s-)dV(s)+\int_{(t_1,t_2]} V(s-)dU(s)+\sum_{u\in(t_1,t_2]} \Delta U_u~\Delta V_u.$$

Introduce cadlag functions of finite variation (on finite intervals):
$$U(s)\defi m([\tau,s])~\mbox{ and } V(s)\defi s-\tau,~~~~~s\in\RR.$$
Then, for $t>0$, applying the previous formulae to $t_1=\tau$ and $t_2=\tau+t$, we see that
\begin{eqnarray}
tm([\tau,\tau+t])&=&\hspace{-1mm}\int\limits_{(\tau,\tau+t]} m([\tau,s))ds+\int\limits_{(\tau,\tau+t]} (s-\tau)~dm(s)=\int\limits_{(0,t]} m([\tau,\tau+s))ds+\int\limits_{(\tau,\tau+t]} (s-\tau)~dm(s)\nonumber\\
&=&\int_{[0,t]} m([\tau,\tau+s))ds+\int_{[\tau,\tau+t]} (s-\tau)~dm(s).\label{e25}
\end{eqnarray}
For the last equality to be proved, it is sufficient to consider a strictly increasing sequence $t_i\uparrow t\ > 0$ and pass to the limit in (\ref{e25}). The case $t\le 0$ is trivial.
\hfill $\Box$

\end{document}